\documentclass[%
aip,
amsmath,amssymb,
preprint,%
]{revtex4-1}
\usepackage{graphicx}% Include figure files
\usepackage{dcolumn}% Align table columns on decimal point
\usepackage{bm}% bold math
\usepackage[utf8]{inputenc}
\usepackage[T1]{fontenc}
\usepackage{mathptmx}
\usepackage{hyperref}
\newcommand{\beq}{\begin{equation}}
\newcommand{\eeq}{\end{equation}}

\begin{document}
	%\preprint{AIP/123-QED}
	\title[]{Effects of forward scattering on the onset of phototactic bioconvection in an algal suspension under diffuse flux without collimated flux}
	%    Force line breaks with \\	
	%\author{S. K. Rajput({\dn s{\rs.} k{\rs.} rAj\7{p}t})}
	\author{S. K. Rajput}
	\altaffiliation[Author to whom correspondence should be addressed. Electronic mail:]{shubh.iiitj@gmail.com.}%Lines break automatically or can be forced with 
	%\author{M. K. Panda({\dn m{\rs.} k{\rs.} p\2A{d}A})}
	\author{M. K. Panda}%
	%\email{Second.Author@institution.edu.}
	\affiliation{Department of Mathematics, PDPM Indian Institute of Information Technology Design and Manufacturing, Jabalpur $482005$, India}%\\This line break forced with \textbackslash\textbackslash	
	%\author{C. Author}
	% \homepage{http://www.Second.institution.edu/~Charlie.Author.}
	%\affiliation{%
	%Second institution and/or address%\\This line break forced% with \\
	%}%	
	%\author{S. K. Rajput}
	% \homepage{http://www.Second.institution.edu/~Charlie.Author.}
	%\affiliation{Department of Mathematics, PDPM Indian Institute of Information Technology Design and Manufacturing, Jabalpur $482005$, India}
	%Second institution and/or address%\\This line break forced% with \\
	%}%	
	%\author{}
	% \homepage{http://www.Second.institution.edu/~Charlie.Author.}
	%\affiliation{Department of Physics, PDPM Indian Institute of Information Technology Design and Manufacturing, Jabalpur $482005$, India}	
	%\author{}
	% \homepage{http://www.Second.institution.edu/~Charlie.Author.}
	%\affiliation{Department of CSE, PDPM Indian Institute of Information Technology Design and Manufacturing, Jabalpur $482005$, India}	
	\date{\today}% It is always \today, today,
	%but any date may be explicitly specified
	%%%%%%%%%%%%%%%%%%%%%%%%%%%%%%%%%%%%%%%%%%%%%%%%%%%%%%	
	\begin{abstract}	
		Phototaxis refers to the directed swimming response influenced by the sensed light intensity in microorganisms. Positive phototaxis involves motion toward the light source, while negative phototaxis entails motion away from it. This study explores the phototactic bioconvection in a suspension of anisotropic scattering phototactic algae, illuminated by diffuse flux without direct collimated flux. The basic state is characterized by zero fluid flow, with the balance between upward and downward swimming due to positive and negative phototaxis, respectively, counteracted by microorganism diffusion. The paper conducts a thorough numerical analysis of linear stability, placing particular emphasis on the impact of forward scattering. The onset of bioconvection manifests either through a stationary mode or an oscillatory mode. The transition between these modes is observed with varying anisotropic coefficients for specific parameter values.		
	\end{abstract}
	%%%%%%%%%%%%%%%%%%%%%%%%%%%%%%%%%%%%%%%%%%%%%%%%%%%%%%	
	\maketitle
	%%%%%%%%%%%%%%%%%%%%%%%%%%%%%%%%%%%%%%%%%%%%%%%%%%%%%%	
	\section{Introduction}
	A remarkable phenomenon where self-propelled microorganisms such as algae, bacteria, protozoa etc. in a fluid exhibit macroscopic convective motion is known as Bioconvection~\cite{20platt1961,23bees2020,24javadi2020}. These microorganisms are slightly denser than the surrounding medium (water) and tend to swim upward on average~\cite{21pedley1992,22hill2005}. However, the biological pattern disappears when the microorganisms stop swimming. On the other hand, up-swimming and greater density are not essential for the development of biological patterns~\cite{21pedley1992}. Instead, swimming microorganisms respond to various environmental stimuli, known as taxes. There are several examples of the taxis such as gravitaxis, chemotaxis, gyrotaxis, and phototaxis etc. Gravitaxis is the swimming response due to gravitational acceleration; chemotaxis is the swimming response due to chemicals; gyrotaxis is caused by the balance between a couple of torques due to gravity and local shear flow; and phototaxis is the swimming response to the light source. This article focuses on the impact of phototaxis on bioconvection.
	
	The wavelength of bioconvection patterns can be influenced by various types of light intensity such as direct (normal and oblique collimated) and indirect (scattered/ diffuse)~\cite{1wager1911,2kitsunezaki2007,3kessler1985,4williams2011,27kage2013,25kessler1986,26kessler1997,28mendelson1998}. When light intensity is high, it can disrupt steady patterns or hinder the formation of patterns in well-stirred cultures~\cite{3kessler1985,4williams2011,5kessler1989}. Changes in the shape and size of bioconvection patterns can also occur due to several reasons. First, photosynthetic motile microorganisms, such as algae, move towards regions of strong light intensity (positive phototaxis) when the light intensity is below a critical value ($G \leq G_c$), and towards regions of low light intensity (negative phototaxis) when it exceeds the critical value ($G \geq G_c$) to avoid photo-damage. This behaviour leads to the accumulation of cells in favourable locations in their natural environment where the light intensity is approximately equal to $G_c$~\cite{6hader1987}. Secondly, the absorption and scattering of light can also impact pattern formation~\cite{7ghorai2010}. Also, diffuse/scattered light can affect pattern formation due to the uniformity of light throughout the suspension.
	
	Phototactic algae can absorb and scatter light and the type of scattering depends on factors such as the size, shape, and refractive index of the cells. Isotropic scattering results in light being scattered uniformly in all directions, while anisotropic scattering is characterized by forward and backward scattering cross-sections. These quantify the amount of light scattered in the forward and backward directions relative to the incident light. Total scattering is the sum of forward and backward scattering. In the case of phototactic algae, the scattered light is mainly in the forward direction due to the size of the cells in relation to the visible wavelength range~\cite{42privoznik1978absorption,43berberoglu2008radiation,13ghorai2013,8panda2020}.
	
	The authors utilized a phototaxis model proposed by Rajput and Panda~\cite{35rajput2023} in this article. This model integrates the Navier–Stokes equations for an incompressible fluid along with a conservation equation for microorganisms and the radiative transfer equation (RTE) for light transport. Scattered/diffuse light is a kind of sunlight created by the scattering of direct (collimated) sunlight by water droplets (clouds) and dust particles present in the atmosphere. During the rainy season, especially on overcast days, particles in clouds disperse sunlight in various directions, resulting in a prevalent diffuse radiation component. The primary form of radiation reaching the Earth's surface becomes diffusive in nature under these conditions. Additionally, regions with high air pollution or atmospheric haze may see the scattering and attenuation of direct sunlight, leading to a higher proportion of diffusive radiation reaching the Earth's surface. Nevertheless, light scatters in a forward direction across an algal suspension, enabling it to penetrate deeper into the suspension. This process has the potential to influence the radiation field, including light intensity and radiative heat flux, by redistributing light intensity profiles. Therefore, this phenomenon governs photosynthesis through phototaxis. It is crucial not to overlook the significance of anisotropic (forward) scattering in understanding phototaxis (photosynthesis) and the resulting bioconvection, especially in the case of dense algal suspensions during the design of efficient photo-bioreactors~\cite{29khan2017,30bees2014}. Many motile algae exhibit strong phototactic behaviour driven by their photosynthetic needs. Therefore, it becomes essential to incorporate the effects of forward scattering into realistic and reliable models of phototaxis to accurately describe the swimming behaviour of such algae. In contrast to Rajput and Panda's work, where the effects of anisotropic (forward) scattering by algae were neglected in their phototaxis model, this study assumes a linearly anisotropic scattering phase function. This assumption acknowledges that phototactic algae scatter light in the forward direction.
	
	In the case of photo-bioconvection in an algal suspension with a finite depth, a sublayer of the suspension is formed where the fluid is stable (motionless), and the movement of cells is influenced by phototaxis and diffusion. The location of this sublayer is determined by the critical intensity of light, $G_c$, and can occur at various depths. If $G<G_c$, the sub-layer forms at the top, and if $G>G_c$, it forms at the bottom. When $G=G_c$ occurs at a certain depth within the suspension, a sub-layer forms in between the suspension, dividing it into two regions. The upper region is gravitationally stable, while the lower region is unstable and experiences fluid flow, which penetrates into the stable upper region. This phenomenon, called penetrative convection, has been observed in various other convection problems~$\cite{9straughan1993,10ghorai2005,11panda2016}$.
	
	\begin{figure}[!htbp]
		\centering
		\includegraphics[width=12cm]{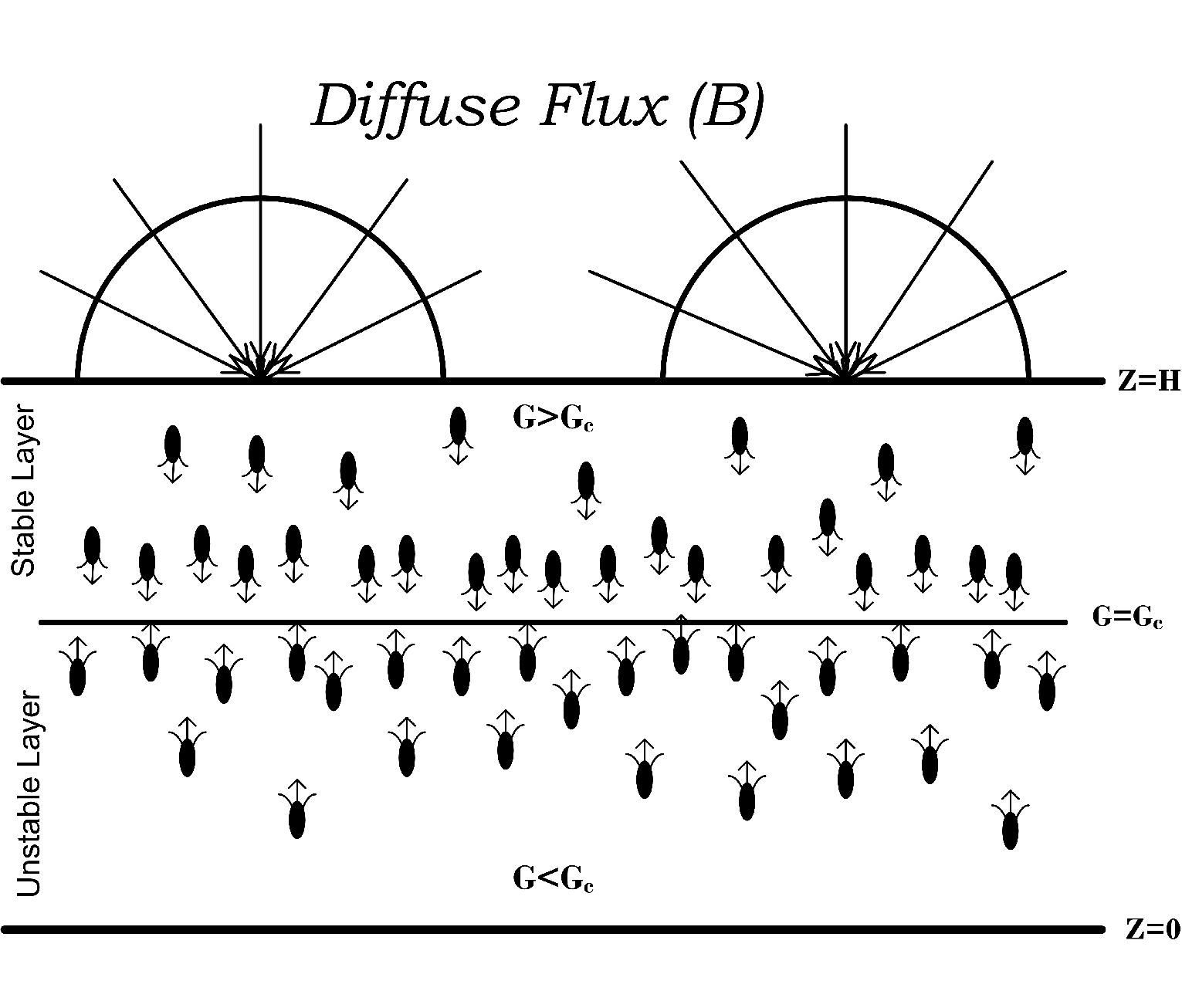}
		\caption{\footnotesize{Formation of sublayer inside the suspension of phototactic algae at a depth where $G=G_c$.}}
		\label{fig1}
	\end{figure}
	
	The phenomenon of phototactic bioconvection in the algal has been extensively explored by researchers under various conditions such as suspension properties and light conditions. The initial investigation by Vincent and Hill~\cite{12vincent1996} focused on the impact of collimated irradiation on an absorbing (non-scattering) medium. Subsequent studies, like that of Ghorai and Hill~\cite{10ghorai2005}, delved into the behaviour of phototactic algal suspensions in two dimensions but did not consider the scattering effect. The effects of light scattering, both isotropic and anisotropic, were investigated by Ghorai $et$ $al$.\cite{7ghorai2010} and Ghorai and Panda\cite{13ghorai2013} under normal collimated irradiation. Panda and Ghorai~\cite{14panda2013} the phenomenon of phototactic bioconvection in algal suspensions 022 5 proposed a model for an isotropically scattering medium in two dimensions, revealing results differing from those of Ghorai and Hill~\cite{10ghorai2005} due to the scattering effect. Panda and Singh~\cite{11panda2016} explored phototactic bioconvection in two dimensions with a non-scattering suspension confined between rigid side walls. The impact of diffuse irradiation, along with collimated irradiation, was studied by Panda $et$ $al$.\cite{15panda2016} in an isotropic scattering medium and by Panda\cite{8panda2020} in an anisotropic medium. In a natural environment, sunlight often strikes the Earth's surface at off-normal angles with oblique irradiation. In this context, the influence of oblique collimated irradiation on non-scattering and isotropic scattering suspensions was examined by Panda $et$ $al$.\cite{16panda2022} and Kumar\cite{17kumar2022}. The effects of diffuse flux alongside oblique collimated flux on isotropic scattering algae suspension were investigated by Panda and Rajput~\cite{35rajput2023}. Recently, Rajput and Panda developed a phototactic model, exploring the effects of scattered/diffuse flux in the absence of direct collimated flux at the onset of bioconvective instability in isotropic scattering algal suspension. However, there has been no study on the onset of phototactic bioconvection in an anisotropic (forward) scattering suspension, incorporating the effects of diffuse/scattered flux in the absence of direct (collimated) irradiation. To gain further insights into phototactic bioconvection, this study investigates the impact of forward (anisotropic) scattering in an algal suspension illuminated by diffuse/scattered flux without direct collimated flux, utilizing the phototaxis model proposed by Rajput and Panda.
	
	The structure of the article consists of several sections. The first section presents the mathematical formulation of the problem, followed by the derivation of a basic solution that represents the equilibrium state. Next, a small perturbation is introduced to the system to study its linear stability by using numerical techniques based on linear perturbation theory. The findings of the model are then described and discussed. Finally, the implications of the model's results are analysed and discussed.
	
	%%%%%%%%%%%%%%%%%%%%%%%%%%%%%%%%%%%%%%%%%%%%%%%%%%%%%% 
	\section{Mathematical formulation of the problem}
	Consider a suspension of forward scattering phototactic algae which is horizontally unbounded and vertically bounded with two parallel horizontal boundaries at $z=0$ and $z=H$. The suspension is assumed to be non-reflective at both the top and bottom boundaries. The top boundary is illuminated uniformly by a diffuse flux, which penetrates throughout the depth of the suspension and illuminates the entire system (as shown in Fig.~\ref{fig1}).
	
	%%%%%%%%%%%%%%%%%%%%%%%%%%%%%%%%%%%%%%%%%%%%%%%%%%%%%%
	\subsection{The mean swimming orientation}
	To calculate the light intensity profile, the Radiative transfer equation (hereafter referred to as RTE) is utilized as a governing equation of light intensity, which is given by
	\begin{equation}\label{1}
	\boldsymbol{s}\cdot\nabla L(\boldsymbol{x},\boldsymbol{s})+(\varkappa+\sigma)L(\boldsymbol{x},\boldsymbol{s})=\frac{\sigma}{4\pi}\int_{4\pi}L(\boldsymbol{x},\boldsymbol{s'})\Xi(\boldsymbol{s},\boldsymbol{s'})d\Omega',
	\end{equation}
	
	where $\varkappa$ and $\sigma$ denote the absorption and scattering coefficients, respectively. The function $\Xi(\boldsymbol{s},\boldsymbol{s'})$ is the scattering phase function. In this study, this function is assumed to be linearly anisotropic. Thus,
	\begin{equation*}
		\Xi(\boldsymbol{s},\boldsymbol{s'})=1+A\cos\theta\cos\theta'
	\end{equation*}
	where $A$ is the linear anisotropic scattering coefficient. Based on the value of this coefficient, we can divide the anisotropic scattering into two parts. First, forward scattering when $0<A\leq 1$, and second, backward scattering when $-1\leq A<0$. If $A=0$, this implies isotropic scattering similar to the previous study.  
	
	%\begin{figure}[!ht]
	%\centering
	%\includegraphics[width=12cm]{diff2.eps}
	%\caption{\footnotesize{Geometric configuration of the problem.}}
	%\label{fig2}
	%\end{figure}
	
	The light intensity at a location $x_b=(x,y,H)$ on the top surface of the suspension is given by
	\begin{equation*}
		L(\boldsymbol{x}_b,\boldsymbol{s})=\frac{B}{\pi}, 
	\end{equation*}
	where $B$ represents the magnitude of diffuse flux (for more detail see Panda $et$ $al$.\cite{15panda2016}). Now, consider  $\varkappa={\varkappa_H}n(\boldsymbol{x})$ and $\sigma={\sigma_s}n(\boldsymbol{x})$, then scattering albedo is defined as $\omega=\sigma_s/\kappa$, where $\kappa=\varkappa_H+\sigma_s$. Therefore, RTE becomes in terms of scattering albedo
	\begin{equation}\label{2}
	\boldsymbol{s}\cdot\nabla L(\boldsymbol{x},\boldsymbol{s})+\kappa n(\boldsymbol{x})L(\boldsymbol{x},\boldsymbol{s})=\frac{\omega \kappa n(\boldsymbol{x})}{4\pi}\int_{4\pi}L(\boldsymbol{x},\boldsymbol{s'})(1+A\cos\theta\cos\theta')d\Omega'.
	\end{equation}
	
	The scattering albedo $\omega\in[0,1]$ and the medium is purely absorbing (or scattering) when $\omega=0$ (or $\omega=1$). 
	
	The total intensity $G(x)$ and radiative heat flux $\boldsymbol{q}(\boldsymbol{x})$ at a fixed point $x$ in the medium is given by
	\begin{equation}
	[G(\boldsymbol{x}),\boldsymbol{q}(\boldsymbol{x})]=\int_{4\pi}L(\boldsymbol{x},\boldsymbol{s})[1,\boldsymbol{s}]d\Omega.
	\end{equation}
	
	The mean cell swimming velocity is defined as
	\begin{equation*}
		\boldsymbol{W}_c=W_c<\boldsymbol{p}>,
	\end{equation*}
	where, $W_c$ is the cell average swimming speed, and $<p>$ represents the average swimming orientation which is given by
	\begin{equation}\label{4}
	<\boldsymbol{p}>=-M(G)\frac{\boldsymbol{{q}}}{|\boldsymbol{{q}}|},
	\end{equation}
	where $M(G)$ is the taxis function such that
	\begin{equation*}
		M(G)=\left\{\begin{array}{ll}\geq 0, & \mbox{if }~~ G(\boldsymbol{x})\leq G_{c},\\< 0, & \mbox{if }~~G(\boldsymbol{x})>G_{c}.  \end{array}\right. 
	\end{equation*}
	
	The specific mathematical expression of $M(G)$ can vary depending on the type of organisms considered~\cite{12vincent1996}.
	
	%%%%%%%%%%%%%%%%%%%%%%%%%%%%%%%%%%%%%%%%%%%%%%%%%%%%%%	
	\subsection{The governing equations}
	Similar to previous continuum models, each cell possesses a volume $v$ and a density of $\rho+\delta\rho$, with $\rho$ representing the density of water. It is assumed that $\delta\rho$ is significantly smaller than $\rho$. The suspension is also considered incompressible. In this context, the equation of continuity can be expressed as follows
	\begin{equation}\label{5}
	\boldsymbol{\nabla}\cdot \boldsymbol{u}=0.
	\end{equation}
	
	The momentum equation
	\begin{equation}\label{6}
	\rho\left(\frac{Du}{Dt}\right)=-\boldsymbol{\nabla} P_e+\mu\nabla^2\boldsymbol{u}-nvg\Delta\rho\hat{\boldsymbol{z}},
	\end{equation}
	where $P_e$ and $\mu$ represent the excess pressure and viscosity of the suspension (equal to the water), respectively.  which is assumed to be equal to the fluid's viscosity.
	
	The governing equation for cell conservation is given by
	\begin{equation}\label{7}
	\frac{\partial n}{\partial t}=-\boldsymbol{\nabla}\cdot{\boldsymbol{J}},
	\end{equation}
	where $\boldsymbol{J}$ is the total cell flux, which is given by 
	\begin{equation}\label{8}
	{\boldsymbol{J}}=n(\boldsymbol{u}+W_c<\boldsymbol{p}>)-\boldsymbol{D}\boldsymbol{\nabla}n.
	\end{equation}
	
	The total cell flux consists of two parts, one due to advection and the other due to diffusion. Here, $\boldsymbol{D}$ is the diffusion tensor, which is considered to be constant and isotropic. Hence, $\boldsymbol{D} = DI$. It is also assumed that the cells exhibit purely phototactic behaviour, so the effect of viscous torque can be neglected. These simplifications allow for the removal of the Fokker-Planck equation from the system of governing equations~\cite{12vincent1996}.
	
	%%%%%%%%%%%%%%%%%%%%%%%%%%%%%%%%%%%%%%%%%%%%%%%%%%%%%%
	\subsection{Boundary conditions}
	In this article, the bottom boundary is considered rigid with no slip while the top boundary is considered stress-free and rigid with no slip. Therefore, the boundary conditions expressed mathematically
	\begin{equation}\label{9}
	\boldsymbol{u}\cdot\hat{\boldsymbol{z}}=\boldsymbol{J}\cdot\hat{\boldsymbol{z}}=0,~~~~~~~~~~z=0,H.
	\end{equation}
	For rigid boundary
	\begin{equation}\label{10}
	\boldsymbol{u}\times\hat{\boldsymbol{z}}=0,~~~~~~~~~~z=0,H,
	\end{equation}
	and for stress-free boundary
	\begin{equation}\label{11}
	\frac{\partial^2}{\partial z^2}(\boldsymbol{u}\cdot\hat{\boldsymbol{z}})=0,~~~~~~~~~~z=H.
	\end{equation}
	
	Now, we consider that the top boundary surface is exposed to a source of diffuse irradiation uniformly, then the boundary condition for intensities are
	\begin{subequations}
		\begin{equation}\label{11a}
		L(x, y, z=H, \theta, \phi) =\frac{B}{\pi},~~~~~ (\pi/2\leq\theta\leq\pi),
		\end{equation}
		\begin{equation}\label{11b}
		L(x, y, z=0, \theta, \phi) =0,~~~~~ (0\leq\theta\leq\pi/2).
		\end{equation}
	\end{subequations}
	
	%%%%%%%%%%%%%%%%%%%%%%%%%%%%%%%%%%%%%%%%%%%%%%%%%%%%%%	
	\subsection{Scaling of the governing equations}
	To scale the governing equations, we use the depth of the layer, $H$, as the length scale, $H^2/D$ as the timescale, $D/H$ as the velocity scale, $H^2/D$ as the pressure scale, and the mean cell concentration, $\bar{n}$, as the concentration scale. After scaling, the governing equations become
	\begin{equation}\label{12}
	\boldsymbol{\nabla}\cdot\boldsymbol{u}=0,
	\end{equation}
	\begin{equation}\label{13}
	S_{c}^{-1}\left(\frac{D\boldsymbol{u}}{Dt}\right)=-\nabla P_e-Rn\hat{\boldsymbol{z}}+\nabla^{2}\boldsymbol{u},
	\end{equation}
	\begin{equation}\label{14}
	\frac{\partial{n}}{\partial{t}}=-{\boldsymbol{\nabla}}\cdot{{\boldsymbol{J}}},
	\end{equation}
	where
	\begin{equation}\label{15}
	{{\boldsymbol{J}}}=n({\boldsymbol{u}}+V_{c}<{\boldsymbol{p}}>)-{\boldsymbol{\nabla}}n.
	\end{equation}
	Here, $S_{c}=\nu/D$ is the Schmidt number, $V_c=W_cH/D$ is the scaled cell swimming speed and $R=\bar{n}vg\Delta{\rho}H^{3}/\nu\rho{D}$ is the Rayleigh number. 
	
	After scaling, the boundary conditions become
	\begin{equation}\label{16}
	\boldsymbol{u}\cdot\hat{\boldsymbol{z}}=\frac{\partial^2}{\partial z^2}(\boldsymbol{u}\cdot\hat{\boldsymbol{z}})=\boldsymbol{J}\cdot\hat{\boldsymbol{z}}=0,~~~~~at~~~~~z=0,
	\end{equation}
	\begin{equation}\label{17}
	\boldsymbol{u}\cdot\hat{\boldsymbol{z}}=u\times\hat{\boldsymbol{z}}=\boldsymbol{J}\cdot\hat{\boldsymbol{z}}=0,~~~~~at~~~~~z=1.
	\end{equation}
	
	In terms of non-dimensional variables, RTE becomes
	\begin{equation}\label{18}
	\boldsymbol{s}\cdot\nabla L(\boldsymbol{x},\boldsymbol{s})+\tau_Hn(\boldsymbol{x})L(\boldsymbol{x},\boldsymbol{s})=\frac{\omega \tau_H n(\boldsymbol{x})}{4\pi}\int_{4\pi}L(\boldsymbol{x},\boldsymbol{s'})(1+A\cos\theta\cos\theta')d\Omega'.
	\end{equation}
	
	where $\tau_H=\kappa\bar{n}H$ denotes the non-dimensional absorption coefficient. In the form of direction cosines $(\xi,\eta,\mu)$ of the direction $\boldsymbol{s}$, where
	\begin{equation*}	\xi=\sin\theta\cos\phi,~~\eta=\sin\theta\sin\phi,~~\mu=\cos\theta
	\end{equation*}
	
	RTE becomes,
	
	\begin{equation}\label{19}
	(\xi,\eta,\mu)\cdot\nabla L(\boldsymbol{x},\boldsymbol{s})+\tau_H nL(\boldsymbol{x},\boldsymbol{s})=\frac{\omega\tau_Hn }{4\pi}\int_{4\pi}L(\boldsymbol{x},\boldsymbol{s'})(1+A\cos\theta\cos\theta')d\Omega',
	\end{equation}
	
	with the boundary conditions
	
	\begin{subequations}
		\begin{equation}\label{20a}
		L(x, y, z=1, \theta, \phi) =\frac{B}{\pi},~~~~~ (\pi/2\leq\theta\leq\pi),
		\end{equation}
		\begin{equation}\label{20b}
		L(x, y, z=0, \theta, \phi) =0,~~~~~ (0\leq\theta\leq\pi/2).
		\end{equation}
	\end{subequations}
	
	%%%%%%%%%%%%%%%%%%%%%%%%%%%%%%%%%%%%%%%%%%%%%%%%%%%%%%	
	\section{The steady solution}
	The steady solution of the system is determined by putting
	\begin{equation*}
		\boldsymbol{u}=0,\quad n=n_s(z),\quad and\quad  L=L_s(z,\theta).
	\end{equation*}
	in the Eqs. $(\ref{15})$-$(\ref{18})$ and $(\ref{25})$.
	Therefore, at the steady state, the total intensity $G_s$ and the radiative heat flux ${\boldsymbol{q}}_s$ can be given by	
	\begin{equation}\label{21}
	G_s=\int_{4\pi}L_s(z,\theta)d\Omega,\quad 
	{{\boldsymbol{q}}}_s=\int_{4\pi}L_s(z,\theta)\boldsymbol{s} d\Omega.
	\end{equation}
	Here, $L_s(z,\theta)$ is independent of $\phi$, so $x$ and $y$ component of $\boldsymbol{q}_s$ vanish. Therefore, $\boldsymbol{q}_s=-q_s\hat{z}$ where $q_s=|\boldsymbol{q}_s|$.
	
	The steady-state intensity $L_s$ can be governed by the equation
	\begin{equation}\label{22}
	\frac{dL_s}{dz}+\left(\frac{\tau_H n_s}{\nu}\right)L_s=\frac{\omega\tau_H n_s}{4\pi\nu}\left(G_s(z)-Aq_s\nu\right).
	\end{equation}
	Here, the light intensity consists solely of a diffuse component denoted as $L_s=L_s^d$, and it is given by
	\begin{equation}\label{23}
	\frac{dL_s^d}{dz}+\left(\frac{\tau_H n_s}{\nu}\right)L_s^d=\frac{\omega\tau_H n_s}{4\pi\nu}\left(G_s(z)-Aq_s\nu\right),
	\end{equation}
	with the boundary conditions
	\begin{subequations}
		\begin{equation}\label{24a}
		L_s^d( z=1, \theta) =\frac{B}{\pi},~~~~~ (\pi/2\leq\theta\leq\pi), 
		\end{equation}
		\begin{equation}\label{24b}
		L_s^d( z=0, \theta) =0,~~~~~(0\leq\theta\leq\pi/2). 
		\end{equation}
	\end{subequations}
	
	Therefore, at the steady state, total intensity $G_s$ can be found by
	
	\begin{equation*}
		G_s=G_s^d=\int_{4\pi}L_s^d(z,\theta)d\Omega.
	\end{equation*}
	
	Now consider a new variable
	
	\begin{equation*}
		\tau(z)=\int_z^1 \tau_H n_s(z')dz',
	\end{equation*}
	
	%*****************************************************
	
	\begin{figure*}
		\includegraphics{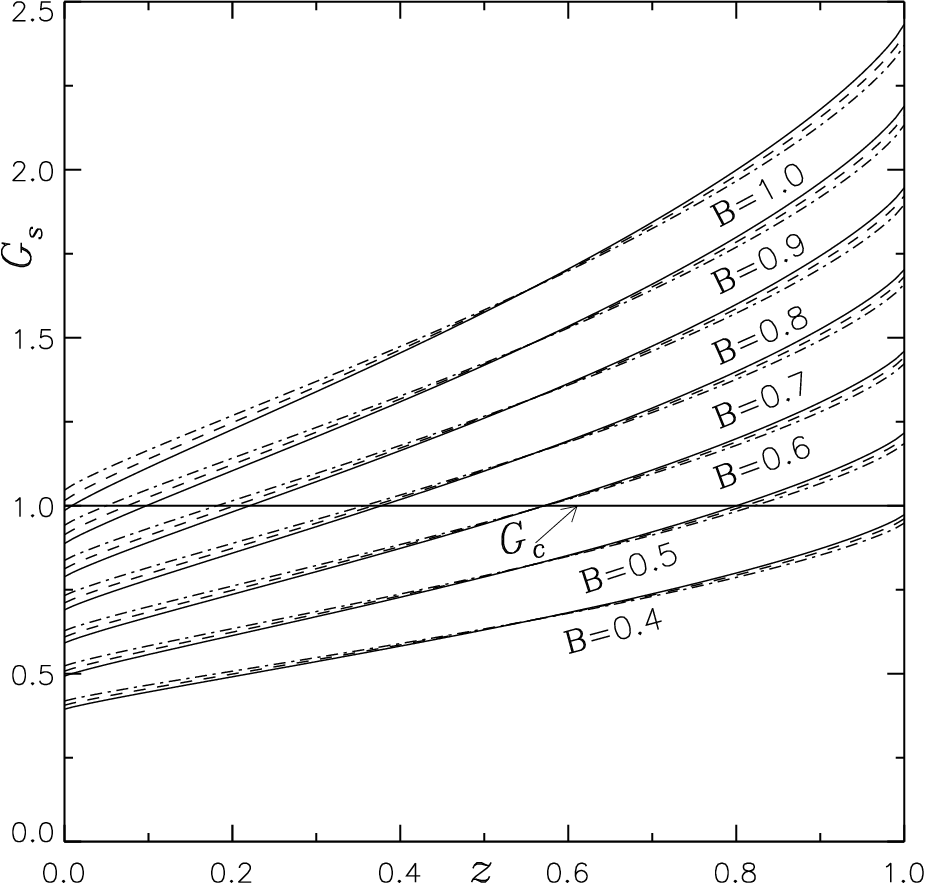}
		\caption{\label{fig3} Exploring the impact of the forward scattering coefficient $A_1$ on the total intensity, $G_s$, within a uniform suspension while maintaining a fixed values $\omega=0.7$ and $\tau_H=0.5$. Various scenarios are considered with different values of B: $A_1=0 (-)$, $A_1=0.4 (---)$, and $A_1=0.8 (-\cdot-\cdot-)$.}
	\end{figure*}
	
	%*****************************************************
	so that, $G_s$ and $\boldsymbol{q}_s$ solely depends on $\tau$ only. Therefore,
	
	\begin{equation}\label{25}
	G_s(\tau)=2BE_2(\tau)+\frac{\omega}{2}\int_0^{\tau_H} \big[E_1(|\tau-\tau'|)G_s(\tau')+A q_s(\tau')E_2(|\tau-\tau'|)sgn(\tau-\tau')\big]d\tau',
	\end{equation}
	\begin{equation}\label{26}
	q_s(\tau)=2BE_3(\tau)+\frac{\omega}{2}\int_0^{\tau_H} \left[A E_3(|\tau-\tau'|)q_s(\tau')+ G_s(\tau')E_2(|\tau-\tau'|)sgn(\tau-\tau')\right]d\tau',
	\end{equation}
	
	where $E_n(x)$ represents the exponential integral of order $n$, and the function $sgn$ stands for the signum function. By using the method of subtraction of singularity~\cite{39crossbie1985}, the couple of Fredholm integral equations (Eq.~(\ref{1}) and Eq.~(\ref{5})) are solved. 
	
	%****************************************************	
	\begin{figure*}
		\includegraphics{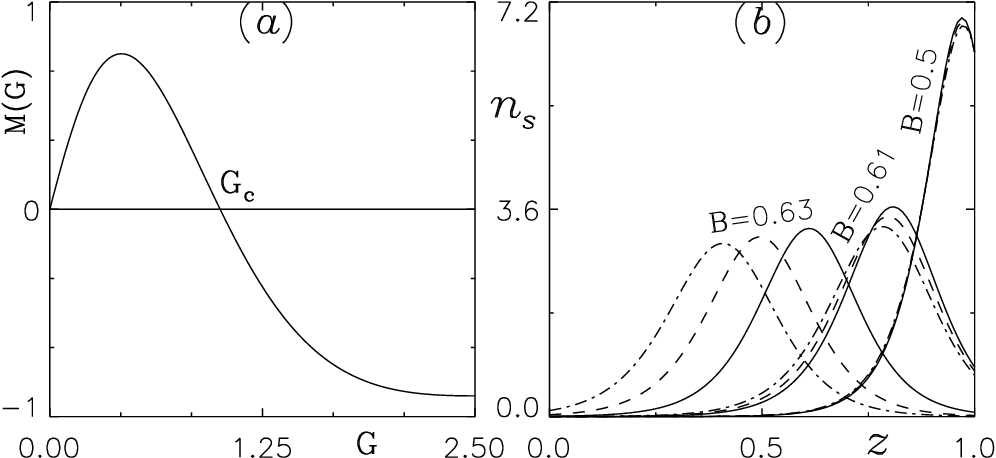}
		\caption{\label{fig4} (a) The photo-response curve featuring a critical intensity value of $G_c=1.0$, and (b) The influence of the forward scattering coefficient, $A_1$, on the fundamental concentration profiles while maintaining a constant $\omega=0.7$, with varied values of B: A$=0 (-)$, A$=0.4 (---)$, and A$=0.8 (-\cdot-\cdot-)$. The parameters $V_c=20$ and $\tau_H=0.5$ are held constant.}
	\end{figure*}
	
	%*****************************************************
	
	In the steady state, the mean swimming orientation can be determined by
	
	\begin{equation*}
		<\boldsymbol{p}_s>=-M(G_s)\frac{{\boldsymbol{q}}_s}{{q}_s}=M_s\hat{z},
	\end{equation*}
	
	where $M_s=M(G_s)$.
	
	The cell concentration, $n_s(z)$, satisfies
	
	\begin{equation}\label{27}
	\frac{dn_s}{dz}-V_cM_sn_s=0,
	\end{equation}
	
	with
	
	\begin{equation}\label{28}
	\int_0^1n_s(z)dz=1.
	\end{equation}
	
	Equations (\ref{25})-(\ref{28}) constitute a boundary value problem that can be solved using a shooting method numerically.
	
	We assume that the microorganisms under consideration are phototactic and similar to Chlamydomonas. Therefore, all parameter values used in this article are the same as used in previous studies (see Table~1 in Panda $et$ $al$.\cite{15panda2016}).
	
	Figure~\ref{fig3} provides insights into how the forward scattering coefficient $A$ affects the total intensity $G_s$ in a uniform suspension for various $B$ values while keeping $\tau_H=0.5$ and $\omega=0.7$ constant. When we increment the value of $A$ from 0 to 0.4 and 0.8, the total intensity $G_s$ changes: it increases at the lower region of the suspension while decreasing at the upper region. Additionally, the position of the critical intensity shifts: it moves downward through the lower half and upward through the upper half of the suspension as $A$ increases.
	
	In Figure~\ref{fig4}, we observe the phototaxis curve for the critical intensity $G_c=1$ and its impact on the steady-state cell concentration, considering parameters $V_c=20$, $\omega=0.7$, and $\tau_H=0.5$ at varying $B$ values. We specifically examine three different $B$ values that lead to maximum concentration occurring at the top, three-quarter height, and mid-height of the suspension in a steady state. Consequently, as we increase the forward scattering coefficient $A$, the location of the maximum concentration in the steady state shifts towards the lower part of the suspension. Furthermore, when the maximum concentration initially occurs at the top of the suspension, its location remains relatively consistent. However, when the maximum concentration shifts towards the three-quarter height and mid-height of the suspension, the distance between these locations increases as $A$ values are raised (see Fig.~\ref{fig4}(b)).
	
	%%%%%%%%%%%%%%%%%%%%%%%%%%%%%%%%%%%%%%%%%%%%%%%%%%%%%%%%%%%%%%%%%%%%%%%%%%%%	
	\section{The linear stability analysis}
	
	For linear stability analysis, we consider a small disturbance 
	(perturbation) to the static state as 
	
	\begin{equation}\label{29}
	\begin{pmatrix}
	\boldsymbol{u}\\
	n\\
	<\boldsymbol{p}>\\
	P_e\\
	L
	\end{pmatrix}
	=
	\begin{pmatrix}
	0\\
	n_s\\
	<\boldsymbol{p}_s>\\
	P_e^s\\
	L_s
	\end{pmatrix}
	+\epsilon
	\begin{pmatrix}
	\boldsymbol{u}_1\\
	n_1\\
	<\boldsymbol{p}_1>\\
	p_e^1\\
	L_1
	\end{pmatrix}
	+\mathcal{O}(\epsilon^2)
	\end{equation}
	
	where, $\boldsymbol{u}_1=(u_1,v_1,w_1)$ and $0<\epsilon\ll 1$.
	
	The perturbed governing equations are given by 
	
	\begin{equation}\label{30}
	\boldsymbol{\nabla}\cdot \boldsymbol{u}_1=0,
	\end{equation}
	
	\begin{equation}\label{31}
	S_{c}^{-1}\left(\frac{\partial \boldsymbol{u_1}}{\partial t}\right)=-\boldsymbol{\nabla} P_e^1-Rn_1\hat{\boldsymbol{z}}+\nabla^{2}\boldsymbol{u}_1,
	\end{equation}
	
	\begin{equation}\label{32}
	\frac{\partial{n_1}}{\partial{t}}+\boldsymbol{\nabla}\cdot(<\boldsymbol{p_s}>n_1+<\boldsymbol{p_1}>n_s)V_c+w_1\frac{dn_s}{dz}=\boldsymbol{\nabla}^2n_1.
	\end{equation}
	
	If $G=G_s+\tilde{\epsilon} G_1+\mathcal{O}(\tilde{\epsilon}^2)=G_s^d+\tilde{\epsilon} G_1^d+\mathcal{O}(\tilde{\epsilon}^2)$, then $G_1$ is given by
	
	\begin{equation}\label{33}
	G_1=G_1^d=\int_{4\pi}L_1^d(\boldsymbol{x},\boldsymbol{s})d\Omega.
	\end{equation}
	
	Similarly
	
	\begin{equation}\label{34}
	{\boldsymbol{q}}_1={\boldsymbol{q}}_1^d=\int_{4\pi}L_1^d(\boldsymbol{x},\boldsymbol{s})\boldsymbol{s}d\Omega.
	\end{equation}
	
	Now, we use the expression
	
	\begin{equation*}
		-M(G_s+{\epsilon} G_1)\frac{{\boldsymbol{q}}_s+{\epsilon}\tilde{\boldsymbol{q}}_1+\mathcal{O}({\epsilon}^2)}{|{\boldsymbol{q}}_s+{\epsilon}\tilde{\boldsymbol{q}}_1+\mathcal{O}({\epsilon}^2)|}-M_s\hat{\boldsymbol{k}},
	\end{equation*}
	
	to find the perturbed swimming direction by collecting $\mathcal{O}({\epsilon} )$ terms. Then, this expression gives
	
	\begin{equation}\label{35}
	<\boldsymbol{p}_1>=G_1\frac{dM_s}{dG}\hat{\boldsymbol{k}}-M_s\frac{\boldsymbol{{q}}_1^{H}}{{q}_s},
	\end{equation}
	
	where ${\boldsymbol{q}}_1^H$ is the horizontal component of perturbed radiative heat flux ${\boldsymbol{q}}_1$.
	
	Now, we substitute the value of $<\boldsymbol{p}_1>$ in equation (\ref{39}) and simplify, obtain
	
	\begin{widetext}
		\begin{equation}\label{36}
		\frac{\partial{n_1}}{\partial{t}}+V_c\frac{\partial}{\partial z}\left(M_sn_1+n_sG_1\frac{dM_s}{dG}\right)-V_cn_s\frac{M_s}{{{q}_s}}\left(\frac{\partial{q}_1^x}{\partial x}+\frac{\partial {q}_1^y}{\partial y}\right)+w_1\frac{dn_s}{dz}=\nabla^2n_1.
		\end{equation}
	\end{widetext}
	
	Now, we eliminate the $P_e^1$ (by taking a double curl of Eq.~\ref{31}) and horizontal component of $\boldsymbol{u}_1$ which reduce the Eqs.~\ref{31}-\ref{32} into two equations for $w_1$ and $n_1$. After that, we account for the normal mode analysis
	
	\begin{equation*}
		\begin{pmatrix}
			w_1\\
			n_1				
		\end{pmatrix}
		=
		\begin{pmatrix}
			W(z)\\
			\Theta(z)
		\end{pmatrix}
		\exp{\big(\sigma t+i(m_1x+m_2y)\big)}.
	\end{equation*}
	
	The perturbed diffuse intensity $L_1^d$ satisfy
	
	\begin{equation}\label{37}
	\zeta\frac{\partial L_1^d}{\partial x}+\mu\frac{\partial L_1^d}{\partial y}+\eta\frac{\partial L_1^d}{\partial z}+\tau_H( n_sL_1^d+n_1L_s^d)=\frac{\omega\tau_H}{4\pi}\bigg(n_sG_1^d+G_s^dn_1+A\eta(n_s\boldsymbol{q}_1\cdot\hat{\boldsymbol{z}}-q_sn_1)\bigg),
	\end{equation}
	
	with the boundary conditions
	
	\begin{subequations}
		\begin{equation}\label{38a}
		L_1^d(x, y, 1, \zeta,\nu, \eta) =0,~~~~(\pi/2\leq\theta\leq\pi), 
		\end{equation}
		\begin{equation}\label{38b}
		L_1^d(x, y, 0, \zeta,\nu, \eta) =0,~~~~ (0\leq\theta\leq\pi/2). 
		\end{equation}
	\end{subequations}
	
	According to the Eq.~(\ref{45}), $L_1^d$ has the following expression
	
	\begin{equation*}
		L_1^d=\Psi^d(z,\zeta,\nu, \eta)\exp{(\sigma t+i(m_1x+m_2y))}. 
	\end{equation*}
	
	From Eq.~(\ref{33}), we get
	
	\begin{widetext}
		\begin{equation*}
			G_1^d=\mathcal{G}^d(z)\exp{(\sigma t+i(m_1x+m_2y))}=\left(\int_{4\pi}\Psi^d(z, \zeta,\nu, \eta)d\Omega\right)\exp{(\sigma t+i(m_1x+m_2y))},
		\end{equation*}
	\end{widetext}
	
	where
	
	\begin{equation}\label{39}
	\mathcal{G}^d(z)=\int_{4\pi}\Psi^d(z, \zeta,\nu, \eta)d\Omega.
	\end{equation}
	
	Similarly, from Eq.(\ref{34}), we get
	
	\begin{equation*}
		\boldsymbol{q}_1=\bigg(q_1^x,q_1^y,q_1^z\bigg)=\bigg(P(z),Q(z),S(z)\bigg)\exp\big(\sigma t+i(m_1x+m_2y)\big)	
	\end{equation*}	
	
	where
	
	\begin{equation}\label{40}
	\bigg(P(z),Q(z),S(z)\bigg)=\int_{4\pi}(\zeta,\nu, \eta)\Psi^d(z, \zeta,\nu, \eta)d\Omega
	\end{equation}
	
	Now $\Psi^d$ satisfies
	
	\begin{widetext}
		\begin{equation}\label{41}
		\frac{d\Psi^d}{dz}+\frac{\big(i(m_1\zeta+m_2\mu)+\tau_H n_s\big)}{\eta}\Psi^d=\frac{\omega\tau_H}{4\pi\eta}\big(n_s\mathcal{G}^d+G_s^d\Theta+A\eta(n_sS-q_s\Theta)\big)-\frac{\tau_H}{\eta}L_s^d\Theta,
		\end{equation}
	\end{widetext}
	
	with the boundary conditions
	
	\begin{subequations}
		\begin{equation}\label{42a}
		\varPsi^d( 1,  \zeta,\nu, \eta) =0,~~~~~ (\pi/2\leq\theta\leq\pi), 
		\end{equation}
		\begin{equation}\label{42b}
		\varPsi^d( 0,\zeta,\nu, \eta) =0,~~~~~ (0\leq\theta\leq\pi/2). 
		\end{equation}
	\end{subequations}
	
	The linear stability equations become
	
	\begin{equation}\label{43}
	\left(\sigma S_c^{-1}+k^2-\frac{d^2}{dz^2}\right)\left( \frac{d^2}{dz^2}-k^2\right)W=Rk^2{\Theta},
	\end{equation}
	
	\begin{widetext}
		\begin{equation}\label{44}
		\left(\sigma+k^2-\frac{d^2}{dz^2}\right){\Theta}+V_c\frac{d}{dz}\left(M_s{\Theta}+n_s\mathcal{G}^d\frac{dM_s}{dG}\right)\\
		-i\frac{V_cn_sM_s}{{q}_s}(m_1P+m_2Q)=-\frac{dn_s}{dz}{W},
		\end{equation}
	\end{widetext}
	
	with the boundary conditions
	
	\begin{equation}\label{45}
	{W}=\frac{d{W}}{dz}=\frac{d{\Theta}}{dz}-V_cM_s{\Theta}-n_sV_c\mathcal{G}^d\frac{dM_s}{dG}=0,~~~~~z=0,1,
	\end{equation}
	
	for stress-free boundary
	
	\begin{equation}\label{46}
	{W}=\frac{d^2{W}}{dz^2}=\frac{d{\Theta}}{dz}-V_cM_s{\Theta}-n_sV_c\mathcal{G}^d\frac{dM_s}{dG}=0,~~~~~z=1.
	\end{equation}
	
	After simplification, Eq.~(\ref{43}) becomes (using D=d/dz)
	
	\begin{equation}\label{47}
	\left(\sigma S_c^{-1}+k^2-D^2\right)\left(D^2-k^2\right)W=Rk^2{\Theta},
	\end{equation}	
	
	\begin{align}\label{48}
		\Upsilon_0(z)+(\sigma+k^2+\Upsilon_1(z)){\Theta}+\Upsilon_2(z)D{\Theta}
		-D^2{\Theta}+Dn_s{W}=0, 
	\end{align}
	
	with the boundary conditions
	
	\begin{equation}\label{49}
	{W}=DW=D\Theta-\Upsilon_2(z){\Theta}-n_sV_c\mathcal{G}^d\frac{dM_s}{dG}=0,~~~~~z=0,1,
	\end{equation}
	
	for stress-free boundary
	
	\begin{equation}\label{50}
	{W}=D^2W=D\Theta-\Upsilon_2(z){\Theta}-n_sV_c\mathcal{G}^d\frac{dM_s}{dG}=0,~~~~~z=1.
	\end{equation}
	
	where
	
	\begin{subequations}
		\begin{equation}\label{51a}
		\Upsilon_0(z)=V_cD\left(n_s\mathcal{G}^d\frac{dM_s}{dG}\right)-i\frac{V_cn_sM_s}{{q}_s}(m_1P+m_2Q),
		\end{equation}
		\begin{equation}\label{51b}
		\Upsilon_1(z)=V_c\frac{dM_s}{dG}DG_s^d,
		\end{equation}
		\begin{equation}\label{51c}
		\Upsilon_2(z)=V_cM_s.
		\end{equation}
	\end{subequations}
	
	%%%%%%%%%%%%%%%%%%%%%%%%%%%%%%%%%%%%%%%%%%%%%%%%%%%%%%
	\section{SOLUTION PROCEDURE}
	
	The NRK method is a numerical technique employed to solve Eqs.~(\ref{47}) and (\ref{48}) iteratively, refining an initial estimate until the desired level of accuracy is reached. This method incorporates specific boundary conditions to calculate the neutral (marginal) stability curves in the $(k, R)$-plane. The neutral curve, denoted as $R^{(n)}(k)$, consists of an infinite number of branches and offers a unique solution to the linear stability problem for a given set of parameters. The most intriguing and unstable bioconvective solution is associated with the branch where the Rayleigh number $R$ reaches its minimum value. Furthermore, the wavelength of the initial disturbance can be determined using the expression $\lambda_c=2\pi/k_c$, where $k_c$ represents the critical wavenumber. Bioconvective solutions are characterized by the formation of vertically aligned convection cells within the suspension. The neutral curve corresponds to points where Re($\gamma$) = 0. If Im($\gamma$) = 0 on this curve, the solution is stationary and non-oscillatory. However, if the imaginary part of $\gamma$ is non-zero, oscillatory solutions may exist. The interplay between stabilizing and destabilizing mechanisms can lead to the emergence of oscillatory solutions. The presence of an overstable solution is identified when the most intriguing solution lies on the oscillatory branch of the marginal stability curve. Assuming the existence of an oscillatory solution, the neutral curve comprises a single oscillatory branch, intersecting the stationary branch at a critical wavenumber of $k=k_b$. This oscillatory branch remains in existence for $k\leq k_b$.
	
	%%%%%%%%%%%%%%%%%%%%%%%%%%%%%%%%%%%%%%%%%%%%%%%%%%%%%    
	\section{NUMERICAL RESULTS}
	
	To investigate how the forward scattering coefficient influences phototactic bioconvection, we conducted a parameter study by varying the value of $A$ within the range of 0 to 0.8. During this investigation, we held several parameters ($S_c$, $V_c$, $\tau_H$, $\omega$, and $B$) constant. To effectively explore a wide range of parameter combinations, we selected a discrete set of constant values to gain insights into the behaviour of phototactic microorganisms. Specifically, we set the Schmidt number as $S_c=20$ and varied the swimming speed of the cells ($V_c$) between 10 and 20, the extinction coefficient ($\tau_H$) from 0.5 to 1.0, the scattering albedo ($\omega$) from 0 to 1, and the intensity of diffuse irradiation ($B$) from 0 to 1. The choice of $B$ was made to ensure that the location of the maximum concentration (referred to as LMC) at the equilibrium state shifted from the top ($z=1$) to the mid-height ($z=1/2$) of the domain, taking into account the effect of the anisotropic scattering coefficient. Within this context, the depth of the region below the LMC that is unstable under gravity (referred to as DGUR) promoted convection, while the height of the gravitationally stable region (referred to as HGSR) above the LMC inhibited fluid motion due to convection. Additionally, the steepness of maximum concentration (referred to as SMC) also played a role in influencing suspension stability. Notably, as DGUR and SMC increased, the critical Rayleigh number ($R_c$) decreased, indicating a decrease in the suspension's stability. Conversely, as HGSR increased, $R_c$ increased, signifying an increase in the suspension's stability.

	\begin{figure*}[!htbp]
		\includegraphics[height=7.2cm,width=16cm]{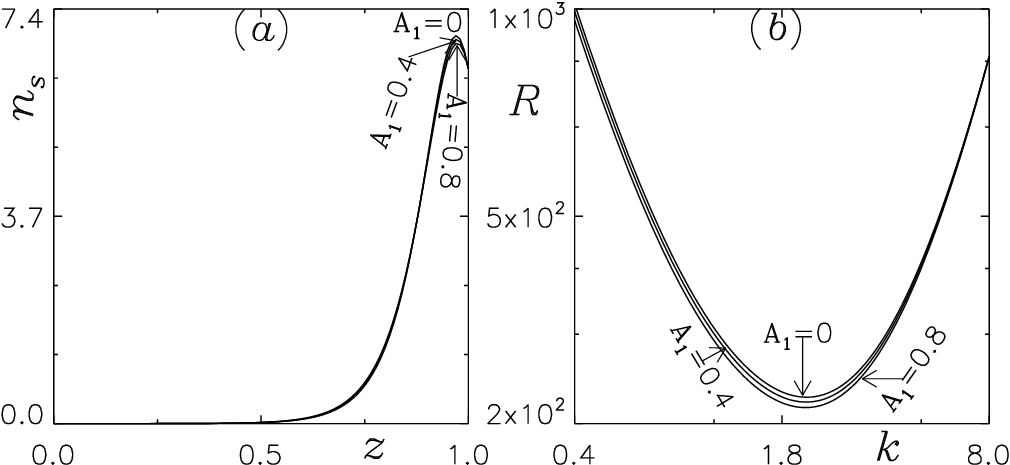}
		\caption{\label{fig5}(a) Examining how the linearly anisotropic scattering coefficient $A_1$, impacts the profiles of base concentration, and (b) exploring the corresponding neutral curves. Here, the parameters that remain fixed are as follows: $V_c = 20$, $\omega = 0.7$, $\tau_H = 0.5$, $B = 0.5$, and $G_c = 1$.}
	\end{figure*}
	
	%%%%%%%%%%%%%%%%%%%%%%%%%%%%%%%%%%%%%%%%%%%%%%%%%%%%%%
	\subsection{When top boundary is stress-free}
	\subsubsection{ $V_c=$ 20 }
	
	(a)~~When $\tau_H=0.5$:
	Figure \ref{fig5} provides an illustration of how the forward scattering coefficient $A_1$ influences the system while maintaining the parameters $V_c=20$, $\tau_H=0.5$, $\omega=0.7$, and $B=0.5$ at fixed values. It offers insights into the base concentration profiles and the associated neutral curves, as seen in Figs. \ref{fig5}(a) and \ref{fig5}(b). In this context, it's important to note that the maximum base concentration is typically found near the top of the suspension. As a consequence, the height of the unstable region (referred to as HUR) remains consistent across all values of $A_1$ since it equals the depth of the suspension itself. Another noteworthy observation is related to the critical depth of the unstable region (known as CDUR), which exhibits minimal variation among different $A_1$ values. However, what does change significantly is the concentration gradient at the upper boundary, and this gradient becomes steeper as $A_1$ is progressively increased. This increase in gradient leads to a reduction in the critical Rayleigh number as $A_1$ is raised from 0 to 0.8, as depicted in Fig. \ref{fig5}.
	
	\begin{figure*}[!htbp]
		\includegraphics[height=7.2cm,width=16cm]{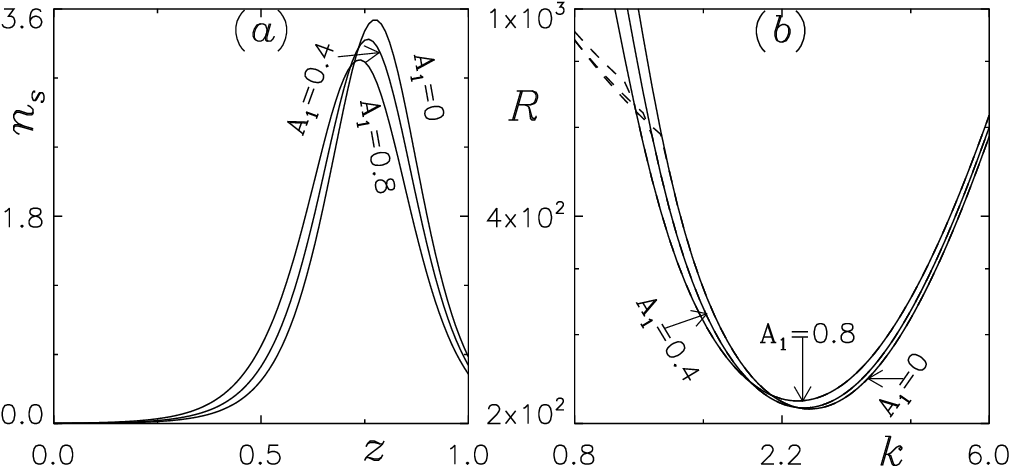}
		\caption{\label{fig6}(a) Examining how the linearly anisotropic scattering coefficient $A_1$, impacts the profiles of base concentration, and (b) exploring the corresponding neutral curves. The parameters that remain fixed are as follows: $V_c = 20$, $\omega = 0.7$, $\tau_H = 0.5$, $B = 0.62$, and $G_c = 1$.}
	\end{figure*}
	
	\begin{figure*}[!bp]
		\includegraphics[height=7.2cm,width=16cm]{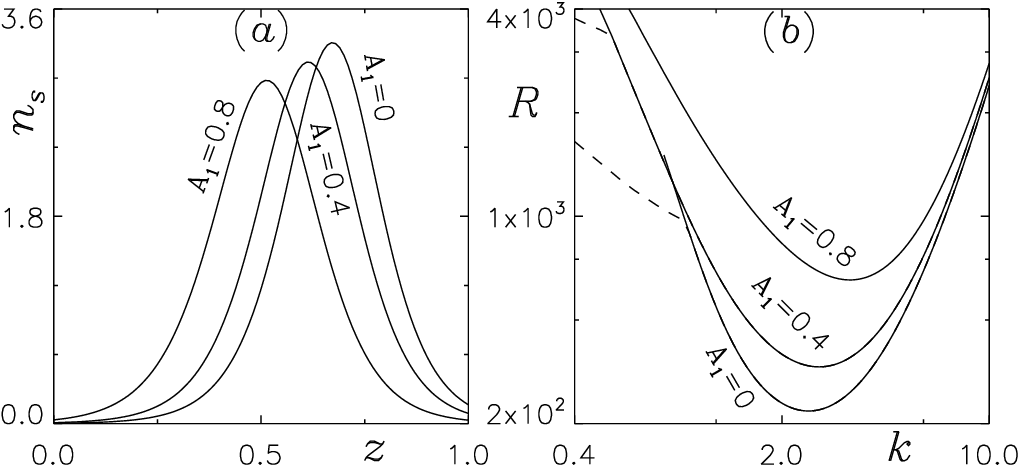}
		\caption{\label{fig7}(a) Examining how the linearly anisotropic scattering coefficient $A_1$, impacts the profiles of base concentration, and (b) exploring the corresponding neutral curves. The parameters that remain fixed are as follows: $V_c = 20$, $\omega = 0.7$, $\tau_H = 0.5$, $B = 0.63$, and $G_c = 1$.}
	\end{figure*}
	
	Figure \ref{fig6} provides insights into how variations in the forward scattering coefficient $A_1$ impact the system's behaviour. This analysis is conducted while keeping the parameters $V_c=20$, $\tau_H=0.5$, $\omega=0.7$, and $B=0.62$ constant. In the results depicted in Fig. \ref{fig6}, we observe that for $A_1=0$, 0.4, and 0.8, oscillatory branches emerge from the stationary branches. However, it's noteworthy that the most unstable solution remains located on the stationary branch for all these values of $A_1$ (as indicated in Fig. \ref{fig6}). This stationary branch corresponds to a base concentration profile where the maximum concentration occurs near the three-quarter height of the suspension. Additionally, as the forward scattering coefficient $A_1$ is increased from 0 to 0.8, both the height of the unstable region (HUR) and the critical depth of the unstable region (CDUR) decrease (see Fig. \ref{fig6}). Furthermore, the steepness of the concentration gradient in the upper stable region exhibits a slight increase as $A_1$ is raised from 0 to 0.8. This change in gradient contributes to the overall increase in the Rayleigh number as the forward scattering coefficient is varied across these values.
	
	Figure \ref{fig7} provides insights into how the forward scattering coefficient $A_1$ impacts the system under specific parameter settings: $V_c=20$, $\tau_H=0.5$, $\omega=0.7$, and $B=0.63$. For this scenario, when $A_1=0$ and 0.4, oscillatory branches emerge from the stationary branches. However, it's essential to note that the most unstable solution remains on the stationary branch for both $A_1=0$ and 0.4 (as indicated in Fig. \ref{fig7}). As $A_1$ increases to 0.8, a stationary solution is observed. In this case, the base concentration profiles are located near the mid-height of the suspension for various $A_1$ values. Furthermore, the height of the unstable region (HUR) and the critical depth of the unstable region (CDUR) both decrease as $A_1$ is increased from 0 to 0.8. However, there is a noticeable increase in the critical Rayleigh number as $A_1$ varies from 0 to 0.8.
	
	\begin{figure*}[!bp]
		\includegraphics[height=7.2cm,width=16cm]{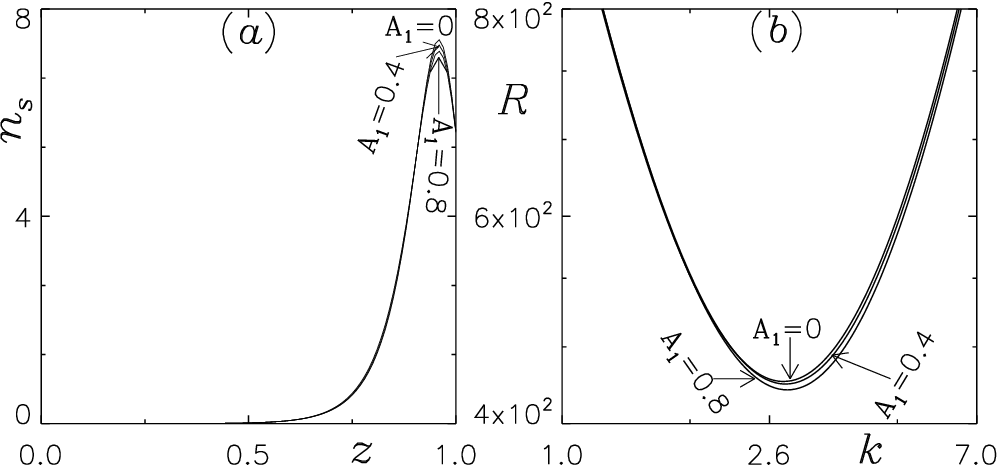}
		\caption{\label{fig8}(a) Examining how the linearly anisotropic scattering coefficient $A_1$, impacts the profiles of base concentration, and (b) exploring the corresponding neutral curves. The parameters that remain fixed are as follows: $V_c = 20$, $\omega = 0.7$, $\tau_H = 1$, $B = 0.6$, and $G_c = 1$.}
	\end{figure*}	
	
	%%%%%%%%%%%%%%%%%%%%%%%%%%%%%%%%%%%%%%%%%%%%%%%%%%%%%%
	
	(b)~~When $\tau_H=1$:	
	In the context of this scenario, where $V_c=20$, $\tau_H=1$, $\omega=0.7$, and $B=0.6$, Figure \ref{fig8} visually demonstrates the impact of the forward scattering coefficient $A_1$ on the base concentration profiles and the associated neutral curves.
	
	\begin{figure*}[!ht]
		\includegraphics[height=7.2cm,width=16cm]{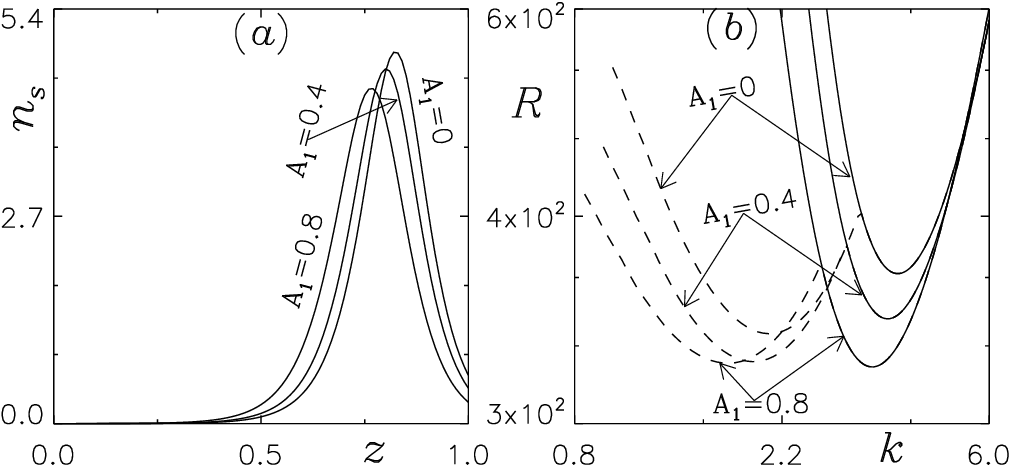}
		\caption{\label{fig9}(a) Examining how the linearly anisotropic scattering coefficient, denoted as $A_1$, impacts the profiles of base concentration, and (b) exploring the corresponding neutral curves. The parameters that remain fixed are as follows: $V_c = 20$, $\omega = 0.7$, $\tau_H = 1$, $B = 0.75$, and $G_c = 1$.}
	\end{figure*}
	
	For all considered values of $A_1$, it is important to note that the height of the unstable region (HUR) is equal to the depth of the suspension. This is because the maximum basic concentration is located near the top. However, as the forward scattering coefficient $A_1$ increases, the concentration gradient at the top of the suspension becomes steeper. Consequently, there is a decrease in the critical Rayleigh number as $A_1$ is raised from 0 to 0.8. It is worth mentioning that the critical depth of the unstable region (CDUR) remains approximately constant across all values of $A_1$. Furthermore, it is essential to emphasize that at the onset of bioconvective instability, the perturbation to the base state is stationary for all values of $A_1$.
	
	\begin{figure*}[!ht]
		\includegraphics[height=15.6cm,width=16cm]{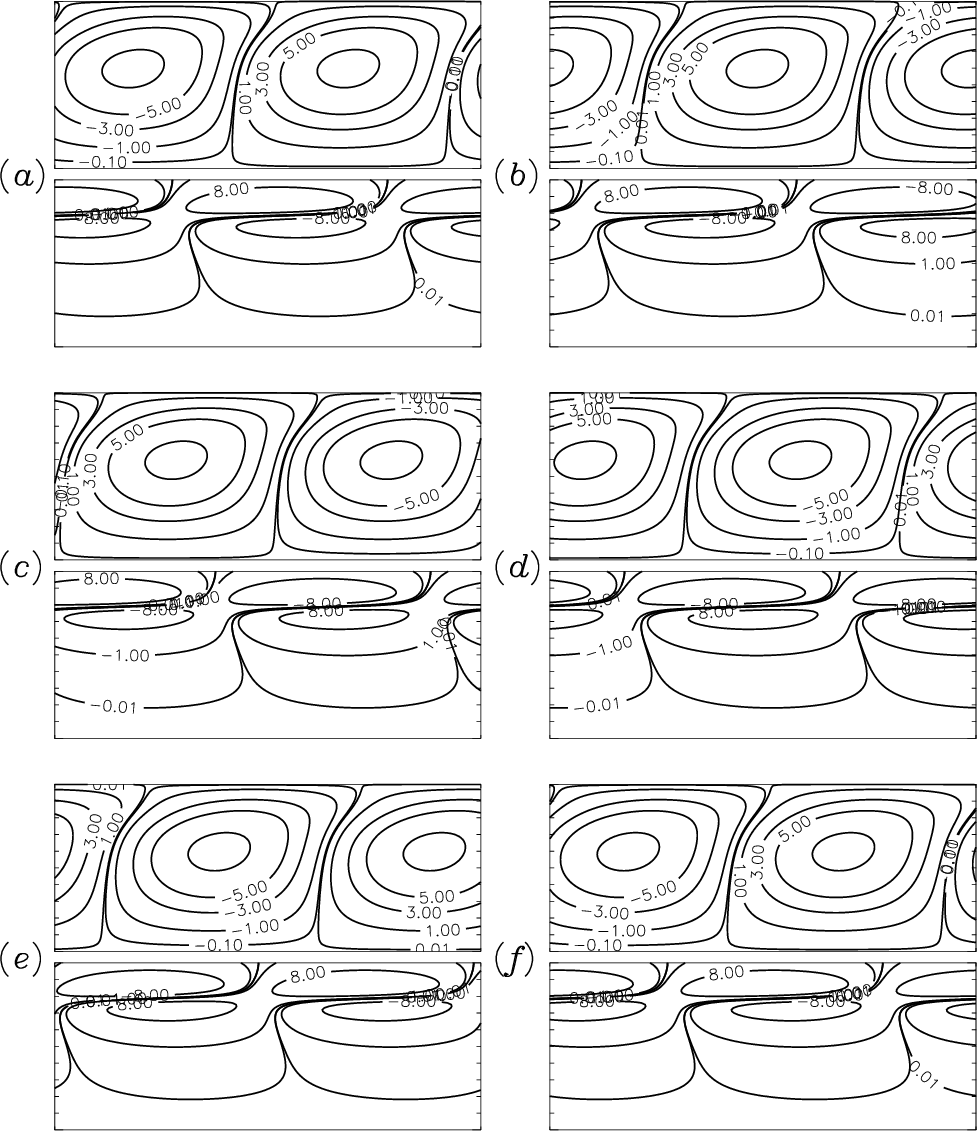}
		\caption{\label{fig10}The flow pattern of perturbed velocity component $w_1$ (top) and perturbed concentration component $n_1$ (bottom) for $A_1=0$ during a cycle of oscillation at initiation of instability. The other parameter values $V_c=20$, $k_H=1$, $\omega=0.7$ and $B=0.75$ are kept constant.}
	\end{figure*}	
	
\begin{figure*}[!htbp]
		\includegraphics[height=6.cm,width=16cm]{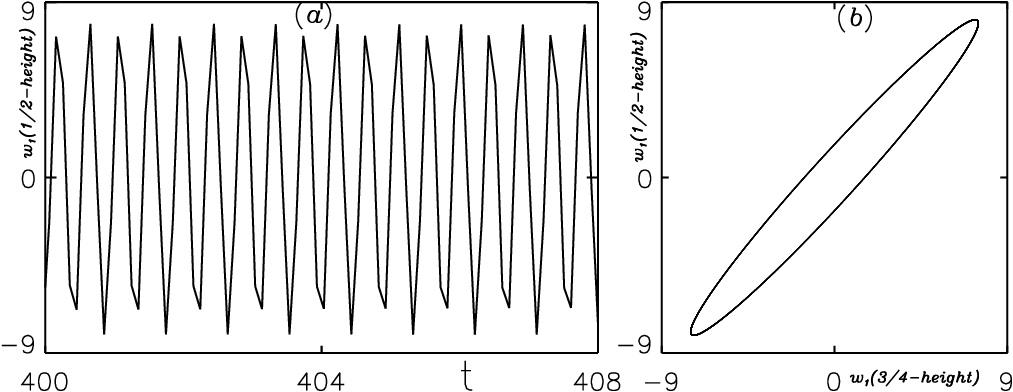}
		\caption{\label{fig11}(a) Variation in the perturbed fluid velocity component $w_1$ with respect to time $t$ and (b) representation of the limit cycle corresponding to the oscillatory branch of the marginal stability curve for a specific forward scattering coefficient $A_1=0$. The governing parameters are maintained at $V_c=20$, $\kappa=1$, $\omega=0.7$, and $B=0.75$.}
	\end{figure*}

	Fig. \ref{fig9} provides insights into how the forward scattering coefficient $A_1$ affects the base concentration profiles and the corresponding neutral curves under constant parameters $V_c=20$, $\tau_H=1$, $\omega=0.7$, and $B=0.75$. In this scenario, oscillatory branches are observed across all values of $A_1$. Specifically, when $A_1=0$, an oscillatory branch emerges from the stationary branch at around $k = k_b \approx 2.75$ and persists for $k \leq k_b$. The most unstable solution, in this case, is located on the oscillatory branch, signifying the onset of overstability at approximately $k_c \approx 1.92$ and $R_c \approx 345$. At this critical point, two eigenvalues are found, forming a complex conjugate pair ($\sigma = 0 \pm 15.78i$), indicating a Hopf (periodic) bifurcation. These eigenvalues correspond to flow patterns that are mirror images of each other and are associated with a wavelength (initial pattern spacing) of approximately $\lambda_c = 3.06$. Notably, nonlinear effects in the bioconvective system become significant on a timescale much shorter than the predicted period of overstability. Consequently, the flow patterns exhibit a travelling wave mechanism, which can be observed through perturbed velocity and concentration components over one cycle of oscillation (see Fig.~\ref{fig10}). Furthermore, this flow destabilization leads to the emergence of a limit cycle due to the period of oscillation, characterized as a Hopf bifurcation with a supercritical nature, ultimately resulting in a stable limit cycle. Fig.~\ref{fig11} shows the predicted time-evolving perturbed ﬂuid velocity component $w_1$ [see Fig. \ref{fig11}(a)] and its corresponding
	phase diagram [see Fig. \ref{fig11}(b)] at $k_c=2.05$.
	
	\begin{figure*}[!bp]
		\includegraphics[height=6.5cm,width=16cm]{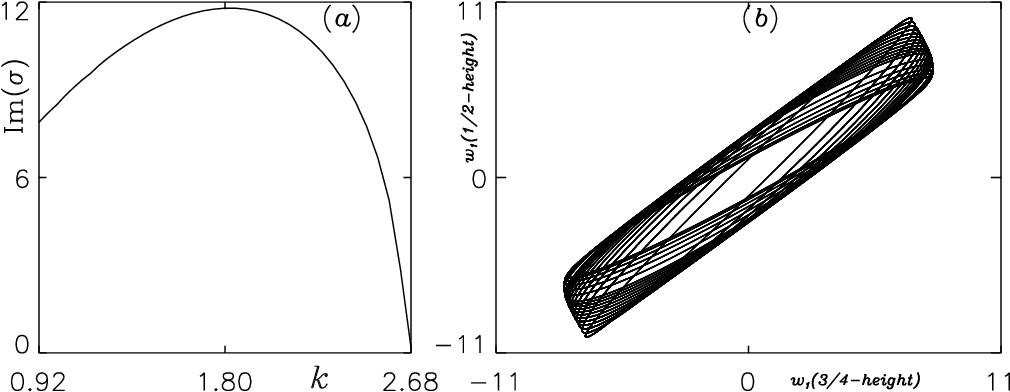}
		\caption{\label{fig12}(a) Variation of $Im(\sigma)$ vs wavenumber $k$ at the initiation of instability, and (b) limit cycle corresponding to the oscillatory branch of marginal stability curve for $A_1=0.8$, where the other governing parameters $V_c=20$, $\tau_H=1$, $\omega=0.7$ and $B=0.75$ are fixed.}
	\end{figure*}	
	
	When $A_1=0.4$, an oscillatory solution is also observed, but for $A_1=0.8$, the critical Rayleigh number lies on the stationary branch (see Fig. \ref{fig9}). Interestingly, both the height of the unstable region (HUR) and the critical depth of the unstable region (CDUR) decreases as $A_1$ increases from 0 to 0.8, promoting convection. However, the steepness of the concentration gradient decreases as well, inhibiting convection. In this case, the inhibiting effect dominates, resulting in an increase in the critical Rayleigh number as $A_1$ varies from 0 to 0.8.
	
	\begin{figure*}[!ht]
		\includegraphics[height=7.2cm,width=16cm]{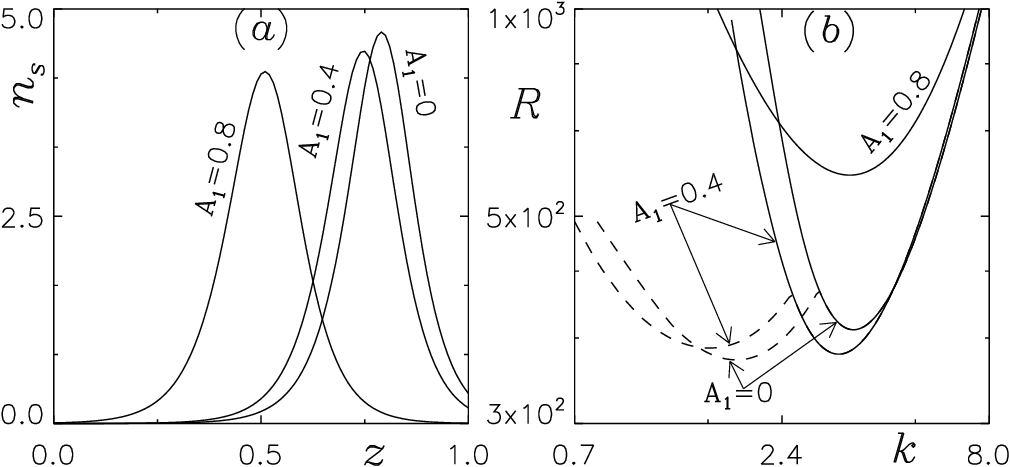}
		\caption{\label{fig13}(a) Examining how the linearly anisotropic scattering coefficient $A_1$, impacts the profiles of base concentration, and (b) exploring the corresponding neutral curves. The parameters that remain fixed are as follows: $V_c = 20$, $\omega = 0.7$, $\tau_H = 1$, $B = 0.76$, and $G_c = 1$.}
	\end{figure*}
	
	\begin{table}[h]
		\caption{\label{tab1} Examining bioconvective solutions while altering the parameter $A_1$ and maintaining constant values for other governing parameters. Results annotated with a double dagger symbol indicate occurrences where the minimum value emerges at the oscillatory branch, whereas a single dagger signifies the presence of oscillations in the $R^{(1)}(k)$ branch of the neutral curve.}
		\begin{ruledtabular}
			\begin{tabular}{ccccccccc}
				$V_c$ & $\tau_H$ & $\omega$ & $B$ & $A_1$ & $\lambda_c$ & $R_c$ & $Im(\sigma)$ & mode\\
				\hline
				20 & 0.5 & 0.7 & 0.5 & 0 & 2.93 & 239.63 & 0 & 1 \\
				20 & 0.5 & 0.7 & 0.5 & 0.4 & 2.93 & 244.27 & 0 & 1\\
				20 & 0.5 & 0.7 & 0.5 & 0.8 & 2.93 & 234.45 & 0 & 1\\
				20 & 0.5 & 0.7 & 0.62 & 0$^{\dagger}$ & 2.52 & 211.67 & 0 & 1\\
				20 & 0.5 & 0.7 & 0.62 & 0.4$^{\dagger}$ & 2.57 & 212.46 & 0 & 1\\
				20 & 0.5 & 0.7 & 0.62 & 0.8$^{\dagger}$ & 2.63 & 218.36 & 0 & 1 \\
				20 & 0.5 & 0.7 & 0.63 & 0$^{\dagger}$ & 2.57 & 272.01 & 0 & 1\\
				20 & 0.5 & 0.7 & 0.63 & 0.4$^{\dagger}$ & 2.35 & 365.36 & 0 & 1\\
				20 & 0.5 & 0.7 & 0.63 & 0.8 & 1.86 & 653.50 & 0 & 2\\
				
				20 & 1 & 0.7 & 0.6 & 0 & 2.23 & 448.41 & 0 &  1\\
				20 & 1 & 0.7 & 0.6 & 0.4 & 2.18 & 446.66 & 0 &  1\\
				20 & 1 & 0.7 & 0.6 & 0.8 & 2.18 & 442.61 & 0 &  1\\
				20 & 1 & 0.7 & 0.75 & 0 & 3.06$^{\ddagger}$ & 362.80$^{\ddagger}$ & 15.78 &  1\\
				20 & 1 & 0.7 & 0.75 & 0.4 & 3.31$^{\ddagger}$ & 345.47$^{\ddagger}$ & 14.05 &  1\\
				20 & 1 & 0.7 & 0.75 & 0.8$^{\dagger}$ & 1.85 & 342.29 & 0 &  1\\
				20 & 1 & 0.7 & 0.76 & 0 & 3.50$^{\ddagger}$ & 332.09$^{\ddagger}$ & 13.31 &  1\\
				20 & 1 & 0.7 & 0.76 & 0.4$^{\dagger}$ & 1.92 & 338.38 & 0 &  1\\
				20 & 1 & 0.7 & 0.76 & 0.8 & 1.80 & 605.70 & 0 &  2\\
			\end{tabular}
		\end{ruledtabular}
	\end{table}
	
	The qualitative behaviour of the bioconvective system is further explored by plotting the phase portrait (bifurcation diagram) of infinitesimal small perturbations. Figure \ref{fig12}(a) illustrates how the positive frequency varies with wavenumber $k$ for the oscillatory branch of the marginal stability curve when $A = 0.8$. The period of oscillation, denoted as $2\pi/Im(\sigma)$, serves as the bifurcation parameter. The bioconvective system exhibits limit cycles or periodic orbits by altering the period of oscillation to higher values through integral multiples of it (see Fig. \ref{fig12}(b)). The phase portrait (bifurcation diagram) includes sequentially decaying orbits, indicating damped oscillations in the bioconvective flow regime as long as the frequency remains non-zero (i.e., when $k < 3.0$) (see Fig. \ref{fig12}(a)). For $k \geq 3.0$, the frequency tends toward zero, and the bioconvective flow regime transitions to stationary (non-oscillatory) convection.
	
	When $V_c=20$, $\tau_H=1$, $\omega=0.7$, and $B=0.76$, we investigate how changes in the forward scattering coefficient $A_1$ impact the base concentration profiles and the neutral curves, as presented in Figure \ref{fig13}. Notably, both the height of the unstable region (HUR) and the critical depth of the unstable region (CDUR) exhibit a decrease as $A_1$ increases from 0 to 0.8. For the specific cases of $A_1=0$ and 0.4, we observe the emergence of oscillatory branches branching out from the stationary branches. However, it's worth noting that the critical Rayleigh number aligns with the oscillatory branch for $A_1=0$ and with the stationary branch for $A_1=0.4$. Consequently, for these values, the system's marginal state assumes an oscillatory behaviour for $A_1=0$ and a stationary one for $A_1=0.4$. Intriguingly, when $A_1$ is further increased to 0.8, a stationary solution becomes prominent. This transition is accompanied by an increase in the critical Rayleigh number as $A_1$ extends from 0 to 0.8 (refer to Fig. \ref{fig13}). The numerical results summarizing the critical Rayleigh number ($R_c$) and the wavelength ($\lambda_c$) for this scenario are presented in Table \ref{tab1}.

	\begin{figure*}[!htbp]
		\includegraphics[height=7.2cm,width=16cm]{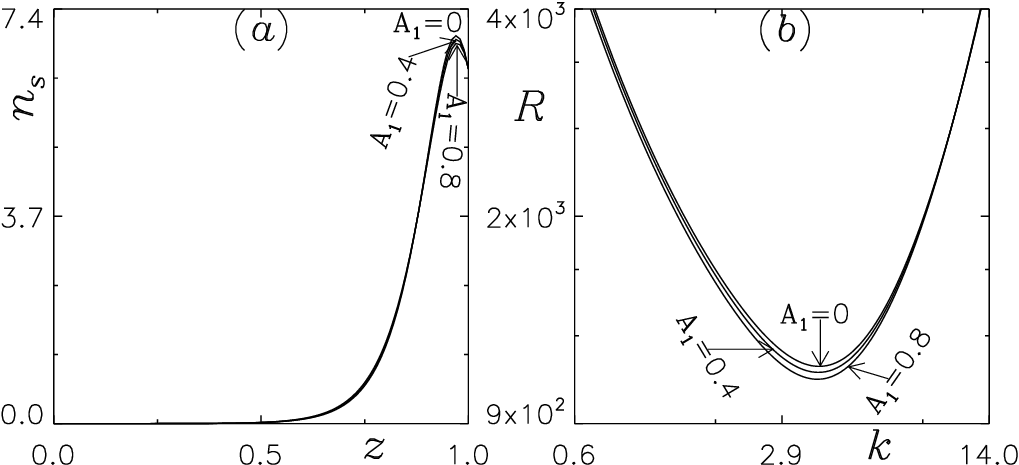}
		\caption{\label{fig14}(a) Examining how the linearly anisotropic scattering coefficient $A_1$, impacts the profiles of base concentration, and (b) exploring the corresponding neutral curves. Here, the parameters that remain fixed are as follows: $V_c = 20$, $\omega = 0.7$, $\tau_H = 0.5$, $B = 0.5$, and $G_c = 1$.}
	\end{figure*}
	\begin{figure*}[!htbp]
		\includegraphics[height=7.2cm,width=16cm]{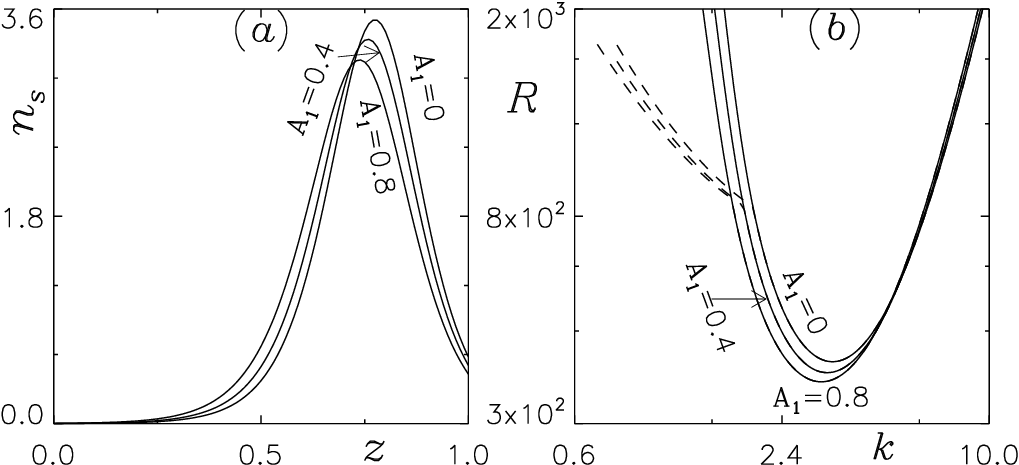}
		\caption{\label{fig15}(a) Examining how the linearly anisotropic scattering coefficient $A_1$, impacts the profiles of base concentration, and (b) exploring the corresponding neutral curves. The parameters that remain fixed are as follows: $V_c = 20$, $\omega = 0.7$, $\tau_H = 0.5$, $B = 0.62$, and $G_c = 1$.}
	\end{figure*}
	
	%%%%%%%%%%%%%%%%%%%%%%%%%%%%%%%%%%%%%%%%%%%%%%%%%%%%%%
	\subsection{When the top boundary is rigid}
	\subsubsection{ $V_c=$ 20 }
	
	(a)~~When $\tau_H=0.5$:
	Figure \ref{fig14} provides an illustration of how the forward scattering coefficient $A_1$ influences the system while maintaining the parameters $V_c=20$, $\tau_H=0.5$, $\omega=0.7$, and $B=0.5$ at fixed values. It offers insights into the base concentration profiles and the associated neutral curves, as seen in Figs. \ref{fig14}(a) and \ref{fig14}(b). In this context, it's important to note that the maximum base concentration is typically found near the top of the suspension. As a consequence, the height of the unstable region (referred to as HUR) remains consistent across all values of $A_1$ since it equals the depth of the suspension itself. Another noteworthy observation is related to the critical depth of the unstable region (known as CDUR), which exhibits minimal variation among different $A_1$ values. However, what does change significantly is the concentration gradient at the upper boundary, and this gradient becomes steeper as $A_1$ is progressively increased. This increase in gradient leads to a reduction in the critical Rayleigh number as $A_1$ is raised from 0 to 0.8, as depicted in Fig. \ref{fig14}.
	
	\begin{figure*}[!bp]
		\includegraphics[height=7.2cm,width=16cm]{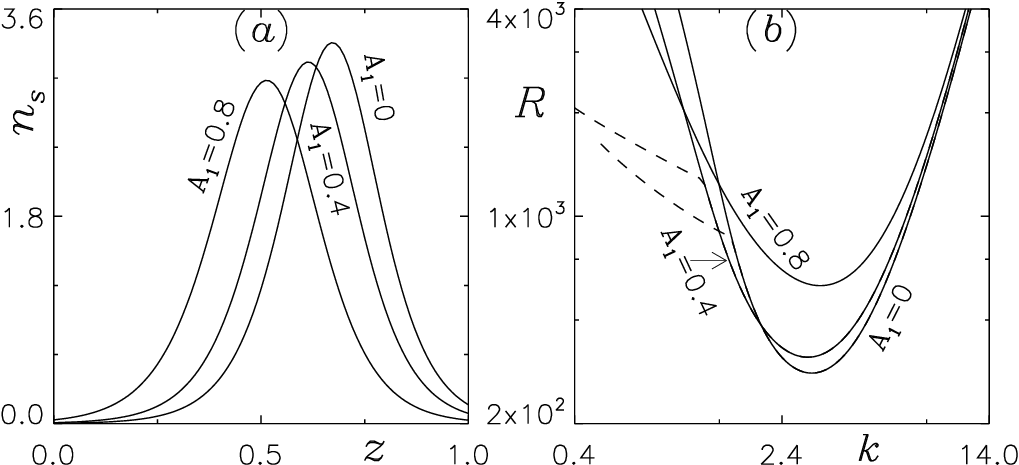}
		\caption{\label{fig16}(a) Examining how the linearly anisotropic scattering coefficient $A_1$, impacts the profiles of base concentration, and (b) exploring the corresponding neutral curves. The parameters that remain fixed are as follows: $V_c = 20$, $\omega = 0.7$, $\tau_H = 0.5$, $B = 0.63$, and $G_c = 1$.}
	\end{figure*}
	
	Figure \ref{fig15} provides insights into how variations in the forward scattering coefficient $A_1$ impact the system's behaviour. This analysis is conducted while keeping the parameters $V_c=20$, $\tau_H=0.5$, $\omega=0.7$, and $B=0.62$ constant. In the results depicted in Fig. \ref{fig15}, we observe that for $A_1=0$, 0.4, and 0.8, oscillatory branches emerge from the stationary branches. However, it's noteworthy that the most unstable solution remains located on the stationary branch for all these values of $A_1$ (as indicated in Fig. \ref{fig15}). This stationary branch corresponds to a base concentration profile, where the maximum concentration occurs near the three-quarter height of the suspension. Additionally, as the forward scattering coefficient $A_1$ is increased from 0 to 0.8, both the height of the unstable region (HUR) and the critical depth of the unstable region (CDUR) decrease (see Fig. \ref{fig15}). Furthermore, the steepness of the concentration gradient in the upper stable region exhibits a slight increase as $A_1$ is raised from 0 to 0.8. This change in gradient contributes to the overall increase in the Rayleigh number as the forward scattering coefficient is varied across these values.
	
	\begin{figure*}[!htbp]
		\includegraphics[height=7.2cm,width=16cm]{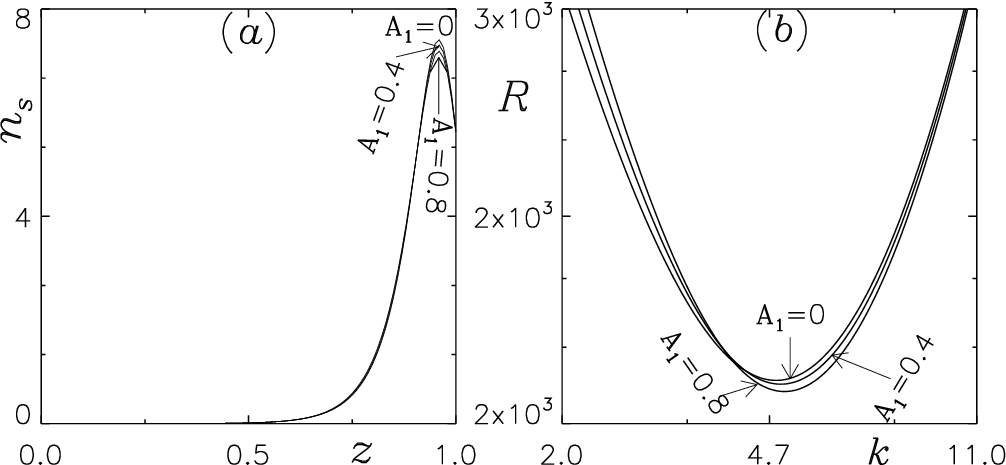}
		\caption{\label{fig17}(a) Examining how the linearly anisotropic scattering coefficient $A_1$, impacts the profiles of base concentration, and (b) exploring the corresponding neutral curves. The parameters that remain fixed are as follows: $V_c = 20$, $\omega = 0.7$, $\tau_H = 1$, $B = 0.6$, and $G_c = 1$.}
	\end{figure*}
	
	Figure \ref{fig16} provides insights into how the forward scattering coefficient $A_1$ impacts the system under specific parameter settings: $V_c=20$, $\tau_H=0.5$, $\omega=0.7$, and $B=0.63$. For this scenario, when $A_1=0$ and 0.4, oscillatory branches emerge from the stationary branches. However, it's essential to note that the most unstable solution remains on the stationary branch for both $A_1=0$ and 0.4 (as indicated in Fig. \ref{fig16}). As $A_1$ increases to 0.8, a stationary solution is observed. In this case, the base concentration profiles are located near the mid-height of the suspension for various $A_1$ values. Furthermore, the height of the unstable region (HUR) and the critical depth of the unstable region (CDUR) both decrease as $A_1$ is increased from 0 to 0.8. However, there is a noticeable increase in the critical Rayleigh number as $A_1$ varies from 0 to 0.8.
	
	\begin{figure*}[!htbp]
		\includegraphics[height=7.2cm,width=16cm]{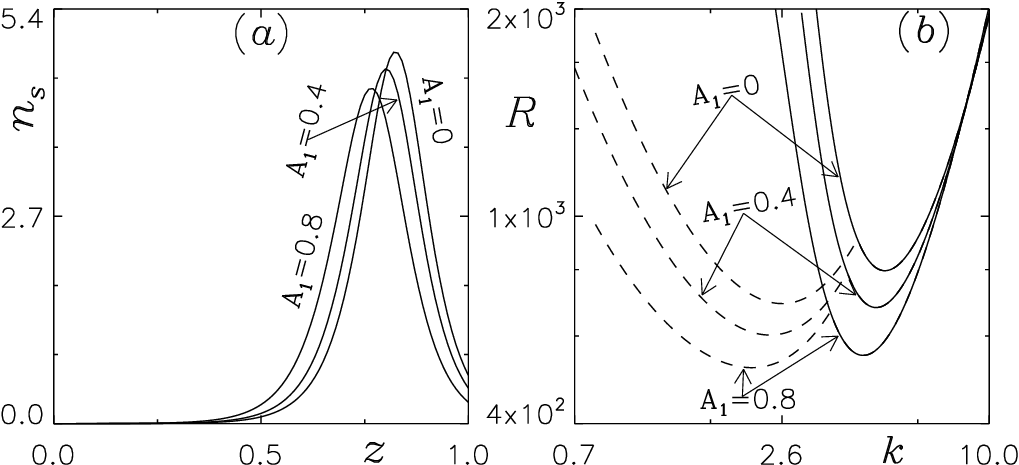}
		\caption{\label{fig18}(a) Examining how the linearly anisotropic scattering coefficient, denoted as $A_1$, impacts the profiles of base concentration, and (b) exploring the corresponding neutral curves. The parameters that remain fixed are as follows: $V_c = 20$, $\omega = 0.7$, $\tau_H = 1$, $B = 0.75$, and $G_c = 1$.}
	\end{figure*}	
	
	%%%%%%%%%%%%%%%%%%%%%%%%%%%%%%%%%%%%%%%%%%%%%%%%%%%%%%
	
	(b)~~When $\tau_H=1$:	
	In the context of this scenario, where $V_c=20$, $\tau_H=1$, $\omega=0.7$, and $B=0.6$, Figure \ref{fig17} visually demonstrates the impact of the forward scattering coefficient $A_1$ on the base concentration profiles and the associated neutral curves.
	
	For all considered values of $A_1$, it is important to note that the height of the unstable region (HUR) is equal to the depth of the suspension. This is because the maximum basic concentration is located near the top. However, as the forward scattering coefficient $A_1$ increases, the concentration gradient at the top of the suspension becomes steeper. Consequently, there is a decrease in the critical Rayleigh number as $A_1$ is raised from 0 to 0.8. It is worth mentioning that the critical depth of the unstable region (CDUR) remains approximately constant across all values of $A_1$. Furthermore, it is essential to emphasize that at the onset of bioconvective instability, the perturbation to the base state is stationary for all values of $A_1$.	
	
	Fig. \ref{fig18} provides insights into how the forward scattering coefficient $A_1$ affects the base concentration profiles and the corresponding neutral curves under constant parameters $V_c=20$, $\tau_H=1$, $\omega=0.7$, and $B=0.75$. In this scenario, oscillatory branches are observed across all values of $A_1$. Specifically, when $A_1=0$, an oscillatory branch emerges from the stationary branch at around $k = k_b \approx 2.75$ and persists for $k \leq k_b$. The most unstable solution, in this case, is located on the oscillatory branch, signifying the onset of overstability at approximately $k_c \approx 1.92$ and $R_c \approx 345$. 
	
	When $A_1=0.4$, an oscillatory solution is also observed, but for $A_1=0.8$, the critical Rayleigh number lies on the stationary branch (see Fig. \ref{fig18}). Interestingly, both the height of the unstable region (HUR) and the critical depth of the unstable region (CDUR) decreases as $A_1$ increases from 0 to 0.8, promoting convection. However, the steepness of the concentration gradient decreases as well, inhibiting convection. In this case, the inhibiting effect dominates, resulting in an increase in the critical Rayleigh number as $A_1$ varies from 0 to 0.8.
	
	\begin{figure*}[!htbp]
		\includegraphics[height=7.2cm,width=16cm]{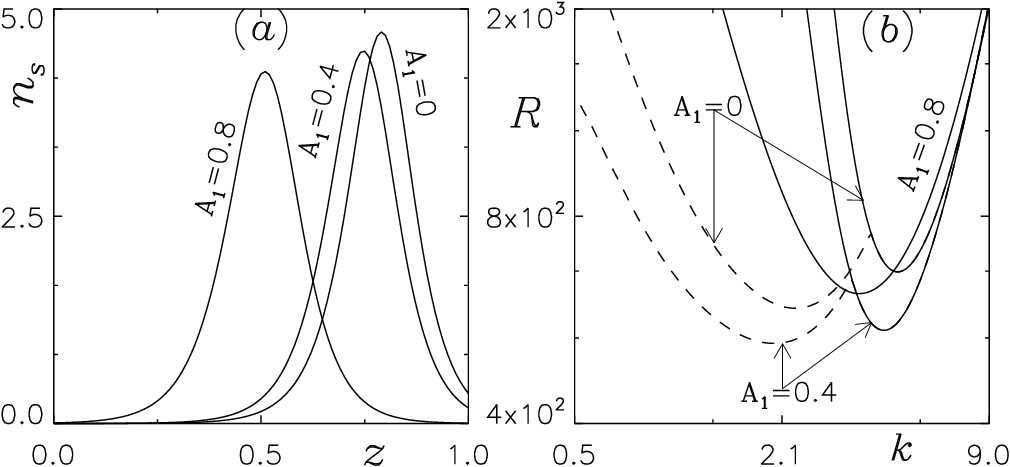}
		\caption{\label{fig22}(a) Examining how the linearly anisotropic scattering coefficient $A_1$, impacts the profiles of base concentration, and (b) exploring the corresponding neutral curves. The parameters that remain fixed are as follows: $V_c = 20$, $\omega = 0.7$, $\tau_H = 1$, $B = 0.76$, and $G_c = 1$.}
	\end{figure*}
	
	\begin{table}[h]
		\caption{\label{tab2} Examining bioconvective solutions while altering the parameter $A_1$ and maintaining constant values for other governing parameters. Results annotated with a double dagger symbol indicate occurrences where the minimum value emerges at the oscillatory branch, whereas a single dagger signifies the presence of oscillations in the $R^{(1)}(k)$ branch of the neutral curve.}
		\begin{ruledtabular}
			\begin{tabular}{ccccccccc}
				$V_c$ & $\tau_H$ & $\omega$ & $B$ & $A_1$ & $\lambda_c$ & $R_c$ & $Im(\sigma)$ & mode\\
				\hline
				20 & 0.5 & 0.7 & 0.5 & 0 & 1.64 & 1111.13 & 0 & 1 \\
				20 & 0.5 & 0.7 & 0.5 & 0.4 & 1.64 & 1086.85 & 0 & 1\\
				20 & 0.5 & 0.7 & 0.5 & 0.8 & 1.66 & 1059.66 & 0 & 1\\
				20 & 0.5 & 0.7 & 0.62 & 0$^{\dagger}$ & 1.83 & 398.22 & 0 & 1\\
				20 & 0.5 & 0.7 & 0.62 & 0.4$^{\dagger}$ & 1.89 & 378.76 & 0 & 1\\
				20 & 0.5 & 0.7 & 0.62 & 0.8$^{\dagger}$ & 1.95 & 363.20 & 0 & 1 \\
				20 & 0.5 & 0.7 & 0.63 & 0$^{\dagger}$ & 2.06 & 350.48 & 0 & 1\\
				20 & 0.5 & 0.7 & 0.63 & 0.4$^{\dagger}$ & 2.14 & 389.80 & 0 & 1\\
				20 & 0.5 & 0.7 & 0.63 & 0.8 & 1.90 & 628.99 & 0 & 2\\
				
				20 & 1 & 0.7 & 0.6 & 0 & 1.31 & 1727.53 & 0 &  1\\
				20 & 1 & 0.7 & 0.6 & 0.4 & 1.28 & 1717.81 & 0 &  1\\
				20 & 1 & 0.7 & 0.6 & 0.8 & 1.26 & 1699.18 & 0 &  1\\
				20 & 1 & 0.7 & 0.75 & 0 & 2.36$^{\ddagger}$ & 69.96$^{\ddagger}$ & 24.71 &  1\\
				20 & 1 & 0.7 & 0.75 & 0.4 & 2.45$^{\ddagger}$ & 621.48$^{\ddagger}$ & 22.25 &  1\\
				20 & 1 & 0.7 & 0.75 & 0.8$^{\dagger}$ & 1.85 & 342.29 & 0 &  1\\
				20 & 1 & 0.7 & 0.76 & 0 & 2.70$^{\ddagger}$ & 581.93$^{\ddagger}$ & 21.43 &  1\\
				20 & 1 & 0.7 & 0.76 & 0.4 & 3.13$^{\ddagger}$ & 519.06$^{\ddagger}$ & 17.51 &  1\\
				20 & 1 & 0.7 & 0.76 & 0.8 & 1.74 & 610.25 & 0 &  2\\
			\end{tabular}
		\end{ruledtabular}
	\end{table}

	When $V_c=20$, $\tau_H=1$, $\omega=0.7$, and $B=0.76$, we investigate how changes in the forward scattering coefficient $A_1$ impact the base concentration profiles and the neutral curves, as presented in Figure \ref{fig22}. Notably, both the height of the unstable region (HUR) and the critical depth of the unstable region (CDUR) exhibit a decrease as $A_1$ increases from 0 to 0.8. For the specific cases of $A_1=0$ and 0.4, we observe the emergence of oscillatory branches branching out from the stationary branches. However, it's worth noting that the critical Rayleigh number aligns with the oscillatory branch for $A_1=0$ and with the stationary branch for $A_1=0.4$. Consequently, for these values, the system's marginal state assumes an oscillatory behaviour for $A_1=0$ and a stationary one for $A_1=0.4$. Intriguingly, when $A_1$ is further increased to 0.8, a stationary solution becomes prominent. This transition is accompanied by an increase in the critical Rayleigh number as $A_1$ extends from 0 to 0.8 (refer to Fig. \ref{fig22}). The numerical results summarizing the critical Rayleigh number ($R_c$) and the wavelength ($\lambda_c$) for this scenario are presented in Table \ref{tab2}.
	
	%%%%%%%%%%%%%%%%%%%%%%%%%%%%%%%%%%%%%%%%%%%%%%%%%%%%%%
	
	\section{CONCLUSIONS}
	
	In conclusion, our investigation has shed light on the impact of forward (anisotropic) scattering in the initiation of phototactic bioconvection within phototactic algae suspensions. This study focused on scenarios with diffuse light as the primary light source, employing a linear anisotropic phase function for scattering.
	
	We observed significant changes when altering the forward scattering coefficient ($A_1$). Increasing $A_1$ led to distinct alterations in the suspension's total intensity profile, with higher values elevating the lower half while reducing the upper half. Moreover, the critical total intensity's position shifted as $A_1$ increased. 
	
	The critical Rayleigh number $R_c$, a key determinant of stability, depends on various parameters, including the unstable region's height, maximum equilibrium concentration, and gradient steepness in the upper stable region. Bioconvection onset exhibited both stationary and oscillatory behaviour, driven by a delicate balance of stabilizing and destabilizing processes, influenced by $B$, concentration gradients, and forward scattering. Forward scattering contributed to convection inhibition by altering the HUR and CDUR. Phototaxis played a dual role in either inhibiting or supporting convection, leading to the emergence of oscillatory solutions.
	
	Notably, under specific parameter conditions, we observed a shift from mode 1 to mode 2 instability in the most unstable solution at bioconvection onset as $A_1$ increased.
	
	The findings also qualitatively align with experimental observations regarding initial pattern wavelengths in the upper half of the suspension, emphasizing the relevance of forward scattering.
	
	While our study provides valuable insights, a comprehensive comparison with quantitative experimental data for pure phototactic algae suspensions is still pending, as most algae species exhibit additional behaviours. Further experimental work is needed to refine our understanding, particularly concerning phototaxis functions, extinction coefficients, and diffusion coefficients. Our study's emphasis on forward scattering in algal suspensions offers a more realistic perspective compared to previous isotropic scattering models. Finally, the proposed phototaxis model holds potential for broader applications beyond this study's scope.
	
	%%%%%%%%%%%%%%%%%%%%%%%%%%%%%%%%%%%%%%%%%%%%%%%%%%%	
	\section*{AUTHOR DECLARATIONS}
	\section*{Conflict of Interest}
	The authors have no conflicts to disclose.
	\section*{DATA AVAILABILITY}
	The data that support the findings of this study are available within the article.
	\bibliography{aniso_diffuse}
	
	%%%%%%%%%%%%%%%%%%%%%%%%%%%%%%%%%%%%%%%%%%%%%%%%%%%%%	
\end{document}